\documentclass[thmsa,nontitlepage,12pt]{article}
\usepackage{rotating}
\usepackage{graphicx}
\usepackage{amsfonts}
\usepackage{amssymb}
\usepackage{amsmath}
\usepackage{latexsym}
\usepackage{sw20elba}
\usepackage{colortbl}
\usepackage[table]{xcolor}
\usepackage{multirow}
\newtheorem{theorem}{Theorem}[section]
\newtheorem{assumption}{Assumption}

\newtheorem{proposition}[theorem]{Proposition}
\newtheorem{remark}[theorem]{Remark}

\numberwithin{equation}{section} \baselineskip=150pt
\input{tcilatex}
\begin{document}

\title{Model Selection for Explosive Models\thanks{%
Yubo Tao, School of Economics, Singapore Management University, 90 Stamford
Road, Singapore 178903. Email: yubo.tao.2014@phdecons.smu.edu.sg. Jun Yu,
School of Economics and Lee Kong Chian School of Business, Singapore
Management University, 90 Stamford Road, Singapore 178903. Email:
yujun@smu.edu.sg.}}
\author{Yubo Tao \ \&\ Jun Yu \ \ \ \ \  \\
\emph{Singapore Management University}}
\maketitle

\begin{abstract}
This paper examines the limit properties of information criteria (such as
AIC, BIC, HQIC) for distinguishing between the unit root model and the
various kinds of explosive models. The explosive models include the
local-to-unit-root model, the mildly explosive model and the regular
explosive model. Initial conditions with different order of magnitude are
considered. Both the OLS estimator and the indirect inference estimator are
studied. It is found that BIC\ and HQIC, but not AIC, consistently select
the unit root model when data come from the unit root model. When data come
from the local-to-unit-root model, both BIC and HQIC\ select the wrong model
with probability approaching 1 while AIC has a positive probability of
selecting the right model in the limit. When data come from the regular
explosive model or from the mildly explosive model in the form of $%
1+n^{\alpha }/n$ with $\alpha \in (0,1)$, all three information criteria
consistently select the true model. Indirect inference estimation can
increase or decrease the probability for information criteria to select the
right model asymptotically relative to OLS, depending on the information
criteria and the true model. Simulation results confirm our asymptotic
results in finite sample.

\medskip

\end{abstract}

\section{Introduction}

Information criteria have found a wide range of practical applications in
empirical work. Examples include choosing explanatory variables in
regression models and selecting lag lengths in time series models.
Frequently used information criteria are AIC of Akaike (1969, 1973), BIC of
Schwarz (1978), HQIC of Hannan and Quinn (1979). A major nice feature in
these information criteria is that the penalty term is trivial to compute
and hence the implementation of them is straightforward and can be made
automatic.

With a growing interest in nonstationarity in time series analysis,
researchers have examined the properties of information criteria in the
context of nonstationary models with the unit root behavior. An important
form of nonstationarity in time series involves explosive roots. Recent
global financial crisis has motivated researchers to study explosive
behavior in economic and financial time series; see, for example, Phillips
and Yu (2011), Phillips, Wu and Yu (2011) and Phillips, Shi and Yu (2015a,
b).

In this paper, we study the limit properties of information criteria for
distinguishing between the unit root model and the explosive models. The
information criteria considered in this paper have a general form and
include AIC, BIC and HQIC as the special cases. The impact of the initial
condition on the limit properties is examined by allowing for an initial
condition of three different orders of magnitude. Moreover, both the OLS
estimator and the indirect inference estimator are studied when
investigating the limit properties of information criteria. The motivation
for the use of indirect inference estimator comes from the existence of
finite sample bias in the OLS estimator and the ability that the indirect
inference method can reduce the bias.

It is found that information criteria consistently choose the unit root
model against the explosive alternatives when data comes from the unit root
model. Second, we prove that the probability for information criteria to
correctly select the explosive model models against the unit root model
depends crucially on both the degree of explosiveness and the size of the
penalty term in information criteria. Finally and surprisingly, we show that
indirect inference estimation can increase or decrease the probability for
information criteria to select the right model asymptotically relative to
OLS, depending on the information criteria and the true model.

The rest of this paper is organized as follows. Section 2 introduces the
models and information criteria, and briefly reviews the literature. Section %
\ref{secLR} gives the limit properties of information criteria for
distinguishing models with an explosive root from the unit root model when
the OLS\ estimator is used. Section 4 gives the limit properties of
information criteria when the indirect inference estimator\ is used. Section %
\ref{secSim} provides Monte Carlo evidence to support the theoretical
results. Section \ref{secCon} concludes. All the detailed proofs are
provided in the appendix. To compress notation, we denote $%
\int\nolimits_{0}^{1}BdB$ and $\int\nolimits_{0}^{1}B^{2}$ in short for $%
\int\nolimits_{0}^{1}B(r)dB(r)$ and $\int\nolimits_{0}^{1}B(r)^{2}dr$
respectively throughout the paper, and $\Rightarrow $ denotes weak
convergence.

\section{Models, Information Criteria and A Literature Review}

The model considered in the present paper is of the form: 
\begin{equation}
X_{t}=\rho _{n}X_{t-1}+u_{t},\text{ }t=1,\cdots ,n,  \label{eqUR}
\end{equation}%
where $u_{t}\overset{iid}{\sim }(0,\sigma ^{2})$ and the model is
initialized at $t=0$ with some $X_{0}$. The autoregressive (AR) coefficient $%
\rho _{n}$ is the crucial parameter that determines the dynamic behavior of $%
X_{t}$. When $\rho _{n}=\rho $ and $\left\vert \rho \right\vert <1$, $X_{t}$
is stationary. When $\rho _{n}=1$, $X_{t}$ has a unit root (UR hereafter).
When $\rho _{n}=1-c_{n}/n=1-c/n$ for $c>0$, $X_{t}$ is near-stationary and
has a root that is local-to-unity (LTUS hereafter) (Phillips, 1987b; Chan
and Wei, 1987). When $\rho _{n}=\rho $ and $\left\vert \rho \right\vert >1$, 
$X_{t}$ has an explosive root (EX hereafter). When $\rho
_{n}=1+c_{n}/n=1+c/n $ for $c>0$, $X_{t}$ is near-explosive and also has a
root that is local-to-unity (LTUE hereafter). When $\rho _{n}=1-c_{n}/n$ for 
$c_{n}\rightarrow \infty $ but $c_{n}/n\searrow 0$, the root represents
moderate deviations from unity and $X_{t}$ is near-stationary (Phillips and
Magdalinos, 2007). When $\rho _{n}=1+c_{n}/n$ for $c_{n}\rightarrow \infty $
but $c_{n}/n\searrow 0$, $X_{t}$ is mildly explosive (hereafter ME).

The asymptotic properties of the OLS\ estimator of the AR coefficient in the
stationary AR(1) model is well known. The rate of convergence is $\sqrt{n}$
and the limiting distribution is Gaussian. Phillips (1987a) provided the
limiting theory for the OLS\ estimator in the UR model and the rate of
convergence is $n$. Phillips (1987b) and Chan and Wei (1987) established the
asymptotic theory for the LTUS and LTUE models. The asymptotic theory is
similar to that in the UR model and the rate of convergence is also $n$. In
the cases of UR and LTU, $u_{t}$ can be weakly dependent stationary.
Anderson (1959) studied the limiting distribution of the OLS\ estimator in
the EX model under the condition that $u_{t}\overset{iid}{\sim }\mathcal{N}%
(0,\sigma ^{2})$ and $X_{0}=0$. The limiting distribution is Cauchy and the
rate of convergence is $\rho ^{n}$. However, no invariance principle
applies. Assuming $X_{0}=o_{p}(\sqrt{n/c_{n}})$, Phillips and Magdalinos
(2007) developed the asymptotic theory for the model with $\rho
_{n}=1-c_{n}/n$ for $c_{n}\rightarrow \infty $ but $c_{n}/n\searrow 0$ and
showed that the asymptotic distribution is invariant to the error
distribution. The rate of convergence is $n/\sqrt{c_{n}}$. If $%
c_{n}=n^{\alpha }$ with $\alpha \in (0,1)$, this rate of convergence bridges
that of UR/LTU models and that of the stationary process. Phillips and
Magdalinos (2007) also developed the asymptotic theory for the ME model. The
rate of convergence is $n\rho _{n}^{n}/c_{n}$. The limiting distribution is
Cauchy which is the same as in the EX model. Interestingly, in the ME case,
the asymptotic theory is independent of the initial condition as long as $%
X_{0}=o_{p}(\sqrt{n/c_{n}})$.

It is known that the OLS estimator of $\rho _{n}$ is biased downward when $%
\rho _{n}=1$ or when $\rho _{n}$ is in the vicinity of unity. In this case,
the indirect inference estimation is effective in reducing the bias.
Phillips (2012) derives the asymptotic theory of the indirect inference
estimator when the model is UR or LTU and $u_{t}\overset{iid}{\sim }\mathcal{%
N}(0,\sigma ^{2})$. The rate of convergence remains unchanged while the
limiting distribution is different from that of the OLS estimator.

Information criteria for model selection have been proposed by Akaike (1969,
1973), Schwarz (1978), Hannan and Quinn (1979), among many others. The
general form of these criteria is

\begin{equation*}
IC_{k}=\log \widehat{\sigma }_{k}^{2}+\frac{kp_{n}}{n},
\end{equation*}%
where $k$ is the number of parameters to be estimated, $\widehat{\sigma }%
_{k}^{2}$ is the estimated $\sigma ^{2}$ when $k$ parameters are estimated.
In general, $IC_{k}$ trades off the term that measures the goodness-of-fit
(i.e. $\log \widehat{\sigma }_{k}^{2}$) and the penalty term that measures
the complexity of the model (i.e. $kp_{n}/n$). Coefficient $p_{n}=2,\log
n,2\log \log n$ corresponds to AIC of Akaike (1973), BIC\ of Schwarz (1978)
and HQIC of Hannan and Quinn (1979). Other forms of $p_{n}$ are possible.

In the time series literature, information criteria have been widely used to
select the lag length both in the family of stationary models and in the
family of nonstationary models; see for example, Ng and Perron (1995) and
Ploberger and Phillips (2003). The information criteria can also be used to
evaluate whether $\rho _{n}=1$ (i.e. $k=0$) or $\rho _{n}\neq 1$ (i.e. $k=1$%
) in Model (\ref{eqUR}). For example, Phillips (2008) obtained limit
properties of $IC_{k}$ for distinguishing between the unit root model and
the stationary model. Phillips and Lee (2015) show that BIC can successfully
distinguish the UR model from the ME model. This is a surprising result as
it is well known that BIC cannot consistently distinguish between the UR\
model and the LTU model; see Ploberger and Phillips (2003).

In this paper we focus our attention to distinguishability between the unit
root model and the three explosive models (i.e., LTUE, ME and EX) after the
candidate models are estimated by OLS or by the indirect inference method.
As a result, we make contributions in two strands of literature, explosive
time series and indirect inference.

To visually understand the difference between the UR\ model, the LTU model
and the ME model, we simulate a sample path of different length ($%
n=100,200,500,1000$) with $y_{0}=0$, based on the same realizations of the
error process, iid $\mathcal{N}(0,1)$, from the following four models, $\rho
_{n}=1$ (UR), $\rho _{n}=1+1/n$ (LTUE), $\rho _{n}=1+n^{0.1}/n$ (ME1), and $%
\rho _{n}=1+n^{0.5}/n$ (ME2). Figures 1-3 give the time series plot of UR
against LTU, UR against ME1, UR against ME2, respectively. It can be seen
from Figure 1 that it is very difficult to distinguish between the UR
process and the LTU process, even when the sample size is as large as 1,000.
When the sample size increases, the gap between the UR\ process and the two
ME processes becomes larger and larger, as apparent in Figure 2 and more so
in Figure 3.

\begin{figure}[htp!]
	\centering
	\includegraphics[width=\textwidth]{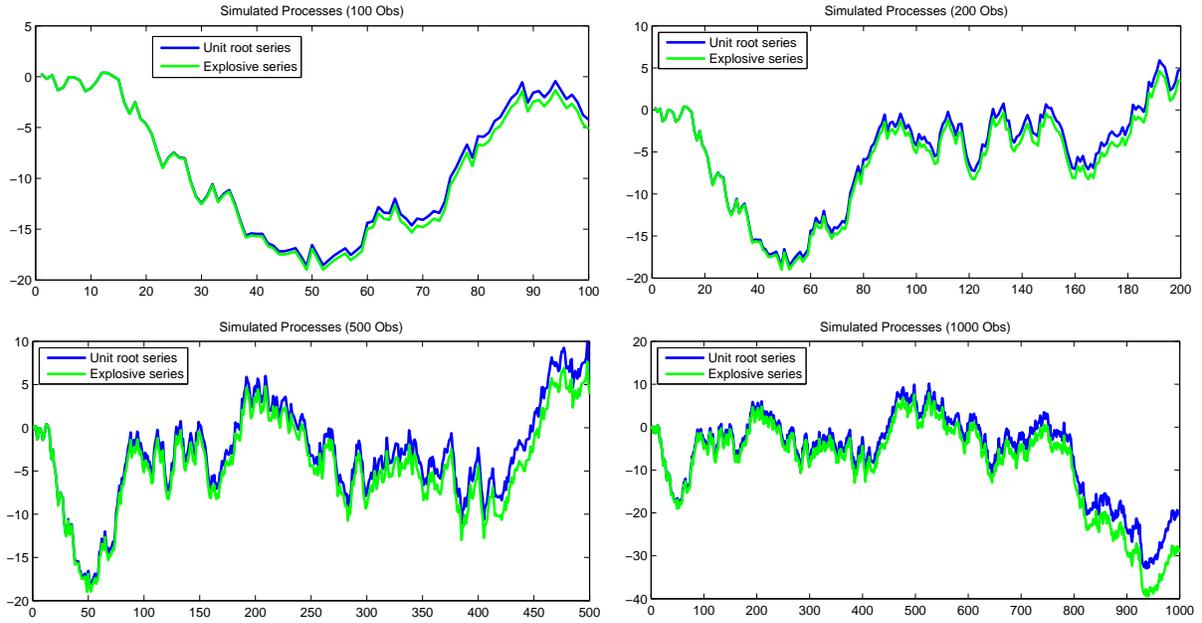}
	\caption{A realization of the UR model
		and the LTU model with $\protect\rho _{n}=1+1/n$.}
\end{figure}

\begin{figure}[htp!]
	\centering
	\includegraphics[width=\textwidth]{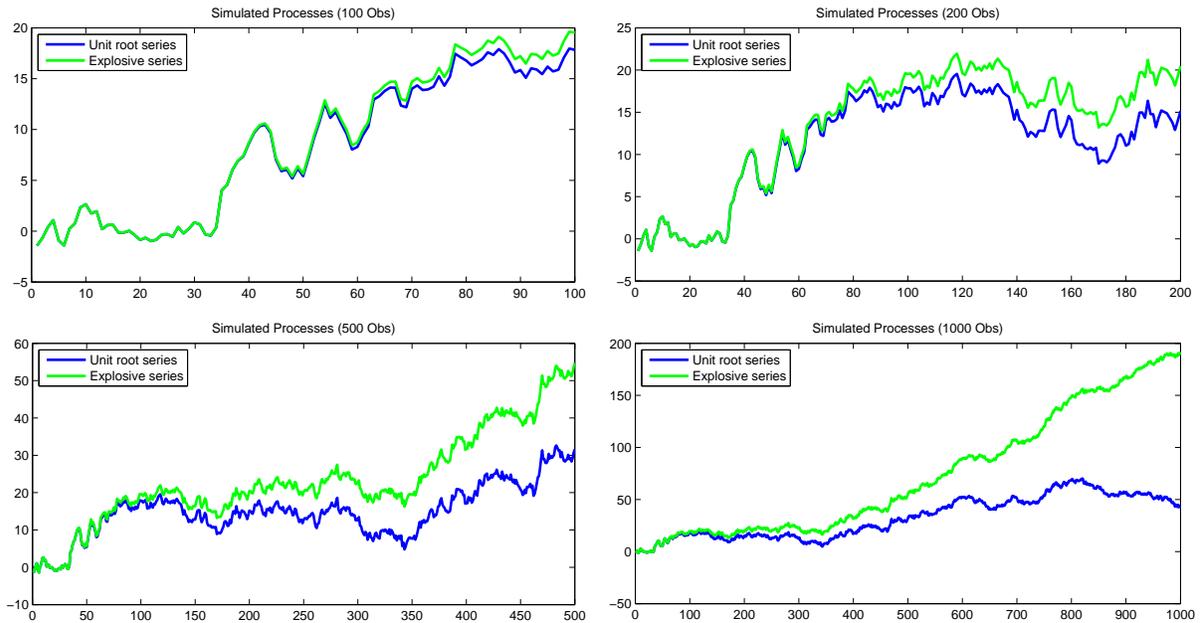}
	\caption{A realization of the UR model
		and the ME process with $\protect\rho _{n}=1+n^{0.1}/n$ (ME1).}
\end{figure}

\begin{figure}[htp!]
	\centering
	\includegraphics[width=\textwidth]{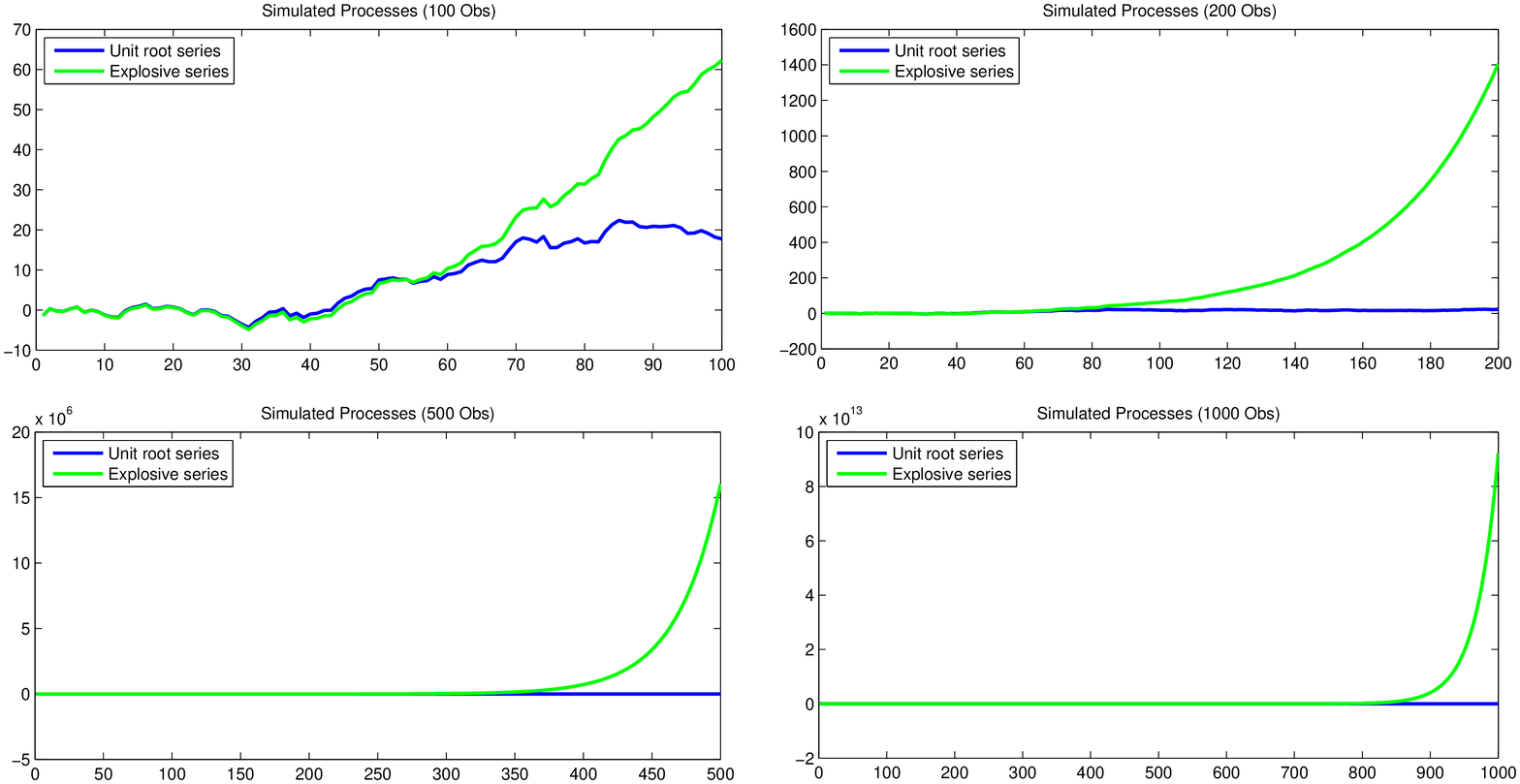}
	\caption{A realization of the UR model
		and the ME model with $\protect\rho _{n}=1+n^{0.5}/n$ (ME2).}
\end{figure}

\section{Limit Properties Based on the OLS Estimator}

\label{secLR}

When the data generating process (DGP) is the UR model, since $\rho _{n}=1$,
we set the parameter count to $k=0$. For the LTU model, the ME model and the
EX model, we need to estimate the AR coefficient and hence set the parameter
count to $k=1$. Throughout the paper we denote $\widehat{\rho }$ the OLS
estimator of $\rho $. $\widehat{k}_{IC}=0$ or $1$ means the information
criterion of the UR\ model is smaller or larger than that of the competing
model when $\rho $ is estimated by OLS. We aim to find the limit of the
following probabilities: 
\begin{align}
& \lim\limits_{n\rightarrow \infty }P\left\{ \widehat{k}_{IC}=0|k=0\right\} ;
\\
& \lim\limits_{n\rightarrow \infty }P\left\{ \widehat{k}_{IC}=1|k=0\right\} ;
\\
& \lim\limits_{n\rightarrow \infty }P\left\{ \widehat{k}_{IC}=0|k=1\right\} ;
\\
& \lim\limits_{n\rightarrow \infty }P\left\{ \widehat{k}_{IC}=1|k=1\right\} .
\end{align}

As shown in Phillips and Magdalinos (2009), the unit root asymptotic
distribution is sensitive to initial conditions in the distant past. To
understand how the initial condition affects the property of $\widehat{k}%
_{IC}$, we follow Phillips and Magdalinos (2009) by assuming alternative
initial conditions.

\begin{assumption}[\textbf{IN}]
\label{assIN}The initial condition has the form 
\begin{equation}
X_{0}(n)=\sum_{j=0}^{\kappa _{n}}u_{-j},
\end{equation}%
where $\kappa _{n}$ is a sequence of integers satisfying $\kappa
_{n}\rightarrow \infty $ and 
\begin{equation}
\frac{\kappa _{n}}{n}\rightarrow \tau \in \left[ 0,\infty \right] \text{, as 
}n\rightarrow \infty .
\end{equation}%
The following cases are distinguished:

\begin{enumerate}
\item[(i)] If $\tau = 0$, $X_0(n) $ is said to be a recent past
initialization.

\item[(ii)] If $\tau \in \left(0, \infty\right)$, $X_0(n) $ is said to be a
distant past initialization.

\item[(iii)] If $\tau =\infty $, $X_{0}(n)$ is said to be an infinite past
initialization.
\end{enumerate}
\end{assumption}

\begin{theorem}
\label{thm01}Under Assumption \ref{assIN} (i) or (ii) or (iii), we have

\begin{enumerate}
\item[(1)] when $p_{n}\rightarrow \infty $ and $p_{n}/n\rightarrow 0$ as $%
n\rightarrow \infty $, 
\begin{align*}
\lim\limits_{n\rightarrow \infty }P\left\{ \widehat{k}_{IC}=0|k=0\right\} &
=\lim\limits_{n\rightarrow \infty }P\left\{ IC_{0}-IC_{1}\leq 0\right\} =1,
\\
\lim\limits_{n\rightarrow \infty }P\left\{ \widehat{k}_{IC}=1|k=0\right\} &
=\lim\limits_{n\rightarrow \infty }P\left\{ IC_{0}-IC_{1}>0\right\} =0.
\end{align*}

\item[(2)] when $p_{n}=2$, the asymptotic distribution under the AIC
criterion is 
\begin{align*}
\lim\limits_{n\rightarrow \infty }P\left\{ \widehat{k}_{AIC}=0|k=0\right\} &
=\lim\limits_{n\rightarrow \infty }P\left\{ AIC_{0}-AIC_{1}\leq 0\right\}
=P\left( \xi ^{2}<2\right) , \\
\lim\limits_{n\rightarrow \infty }P\left\{ \widehat{k}_{AIC}=1|k=0\right\} &
=\lim\limits_{n\rightarrow \infty }P\left\{ AIC_{0}-AIC_{1}>0\right\}
=1-P\left( \xi ^{2}<2\right) .
\end{align*}%
where 
\begin{equation*}
\xi ^{2}=%
\begin{cases}
\dfrac{\left( \int_{0}^{1}BdB\right) ^{2}}{\int_{0}^{1}B^{2}}, & \text{if }%
\tau =0 \\ 
\dfrac{\left( \int_{0}^{1}B_{\tau }dB\right) ^{2}}{\int_{0}^{1}B_{\tau }^{2}}%
, & \text{if }\tau \in (0,\infty ) \\ 
B(1)^{2}, & \text{if }\tau =\infty%
\end{cases}%
,
\end{equation*}%
with $B(s)$ being a Brownian motion, and 
\begin{equation*}
B_{\tau }(s)=B(s)+\sqrt{\tau }B_{0}(1),
\end{equation*}%
with $B_{0}(s)$ being an independent Brownian motion.
\end{enumerate}
\end{theorem}

\begin{remark}
\label{rmk01}Theorem \ref{thm01} is the same as Theorem 1 in Phillips (2008)
for distinguishing between the UR model and the stationary model. The
condition that $p_{n}\rightarrow \infty $ and $p_{n}/n\rightarrow 0$ covers
BIC and HQIC and hence, both BIC and HQIC can consistently select the UR
model. The AIC criterion is inconsistent and its asymptotic distribution
depends on $\xi ^{2}$, the squared unit root $t$-statistic for the OLS
estimator.
\end{remark}

\begin{remark}
\label{rmk01b}The validity of Theorem \ref{thm01} does not require the iid
assumption for the error term $u_{t}$. If we follow Phillips (2008) by
denoting $F(L)=\sum_{j=0}^{\infty }F_{j}L^{j}$, with $F_{0}=1$ and $F(1)\neq
0$, and letting $u_{s}$ have Wold representation 
\begin{equation}
u_{s}=F(L)\varepsilon _{s}=\sum_{j=0}^{\infty }F_{j}\varepsilon _{s-j}\text{%
, with }\sum_{j=0}^{\infty }j^{1/2}\left\vert F_{j}\right\vert <\infty ,
\label{dep}
\end{equation}%
where $\varepsilon _{t}\overset{iid}{\sim }\left( 0,\sigma _{\varepsilon
}^{2}\right) $, the results in Theorem \ref{thm01} continue to hold.
However, both $B_{0}$ and $\xi ^{2}$ need to be modified to accommodate the
dependence in $u_{t}$ as in Phillips (2008).
\end{remark}

\begin{theorem}
\label{thm02a}Let Assumption \ref{assIN} (i) or (ii) holds. Assume the true
DGP is the LTUE model.

\begin{enumerate}
\item[(1)] When $p_{n}\rightarrow \infty $ and $p_{n}/n\rightarrow 0$ as $%
n\rightarrow \infty $, 
\begin{align*}
\lim\limits_{n\rightarrow \infty }P\left\{ \widehat{k}_{IC}=0|k=1\right\} &
=\lim\limits_{n\rightarrow \infty }P\left\{ \dfrac{n}{p_{n}}\left(
IC_{1}-IC_{0}\right) >0\right\} =1, \\
\lim\limits_{n\rightarrow \infty }P\left\{ \widehat{k}_{IC}=1|k=1\right\} &
=\lim\limits_{n\rightarrow \infty }P\left\{ \dfrac{n}{p_{n}}\left(
IC_{1}-IC_{0}\right) \leq 0\right\} =0.
\end{align*}

\item[(2)] When $p_{n}=2$, the asymptotic distribution of the AIC criterion
is 
\begin{align*}
\lim\limits_{n\rightarrow \infty }P\left\{ \widehat{k}_{AIC}=0|k=1\right\} &
=\lim\limits_{n\rightarrow \infty }P\left\{ {n}\left( AIC_{1}-AIC_{0}\right)
>0\right\} =1-P\left( \zeta ^{2}>2\right) , \\
\lim\limits_{n\rightarrow \infty }P\left\{ \widehat{k}_{AIC}=1|k=1\right\} &
=\lim\limits_{n\rightarrow \infty }P\left\{ {n}\left( AIC_{1}-AIC_{0}\right)
\leq 0\right\} =P\left( \zeta ^{2}>2\right) ,
\end{align*}%
where 
\begin{equation*}
\zeta ^{2}=\frac{\left( \int_{0}^{1}J_{c}dB\right) ^{2}}{%
\int_{0}^{1}J_{c}^{2}}+2{c}\int_{0}^{1}J_{c}dB+c^{2}\int_{0}^{1}J_{c}^{2},
\end{equation*}%
with 
\begin{equation*}
J_{c}(r)=\int_{0}^{r}\exp \left\{ c(r-s)\right\} dB(s).
\end{equation*}
\end{enumerate}
\end{theorem}

\begin{remark}
\label{rmk02a}Theorem \ref{thm02a} shows that all the information criteria
are inconsistent in distinguishing between the LTUE model and the UR models
when data comes from the LTUE model. AIC selects the wrong model with
probability going to $1-P\left( \zeta ^{2}>2\right) $, which depends on the
localization constant $c$. This problem worsens for BIC and HQIC as the
probability of selecting the wrong model goes to one. Note that BIC is well
known to be blind to local alternatives; see, for example, Ploberger and
Phillips (2003).
\end{remark}

\begin{theorem}
\label{thm02b}Let Assumption \ref{assIN} (i) or (ii) holds. Assume the true
DGP is the ME model.

\begin{enumerate}
\item[(1)] When $\lim\limits_{n\rightarrow \infty }\dfrac{p_{n}}{\rho
_{n}^{2n}}=0,$ 
\begin{align*}
\lim\limits_{n\rightarrow \infty }P\left\{ \widehat{k}_{IC}=0|k=1\right\} &
=\lim\limits_{n\rightarrow \infty }P\left\{ \dfrac{n}{\rho _{n}^{2n}}\left(
IC_{1}-IC_{0}\right) >0\right\} =0, \\
\lim\limits_{n\rightarrow \infty }P\left\{ \widehat{k}_{IC}=1|k=1\right\} &
=\lim\limits_{n\rightarrow \infty }P\left\{ \dfrac{n}{\rho _{n}^{2n}}\left(
IC_{1}-IC_{0}\right) \leq 0\right\} =1.
\end{align*}

\item[(2)] When $\lim\limits_{n\rightarrow \infty }\dfrac{p_{n}}{\rho
_{n}^{2n}}=\pi \in (0,+\infty ),$ 
\begin{align*}
\lim\limits_{n\rightarrow \infty }P\left\{ \widehat{k}_{IC}=0|k=1\right\} &
=\lim\limits_{n\rightarrow \infty }P\left\{ \dfrac{n}{\rho _{n}^{2n}}\left(
IC_{1}-IC_{0}\right) >0\right\} =P\left( \chi ^{2}(1)<4\pi \right) , \\
\lim\limits_{n\rightarrow \infty }P\left\{ \widehat{k}_{IC}=1|k=1\right\} &
=\lim\limits_{n\rightarrow \infty }P\left\{ \dfrac{n}{\rho _{n}^{2n}}\left(
IC_{1}-IC_{0}\right) \leq 0\right\} =1-P\left( \chi ^{2}(1)<4\pi \right) .
\end{align*}

\item[(3)] When $\lim\limits_{n\rightarrow \infty }\dfrac{p_{n}}{\rho
_{n}^{2n}}\rightarrow +\infty ,$ 
\begin{align*}
\lim\limits_{n\rightarrow \infty }P\left\{ \widehat{k}_{IC}=0|k=1\right\} &
=\lim\limits_{n\rightarrow \infty }P\left\{ \dfrac{n}{p_{n}}\left(
IC_{1}-IC_{0}\right) >0\right\} =1, \\
\lim\limits_{n\rightarrow \infty }P\left\{ \widehat{k}_{IC}=1|k=1\right\} &
=\lim\limits_{n\rightarrow \infty }P\left\{ \dfrac{n}{p_{n}}\left(
IC_{1}-IC_{0}\right) \leq 0\right\} =0.
\end{align*}
\end{enumerate}
\end{theorem}

\begin{remark}
\label{rmk02b}Theorem \ref{thm02b} shows that the limit probability of
selecting the correct model by information criteria under the ME model
depends critically on two parameters, $c_{n}$, $p_{n}$. As expected, the
larger $c_{n}$, the further the model away from the UR\ model and the higher
probability for the information criteria to select the correct model.
Interestingly, the smaller $p_{n}$, the higher probability for the
information criteria to select the correct model. From Phillips and
Magdalinos (2009), we know $\rho _{n}^{-n}=o(c_{n}^{-1})$ and hence $\rho
_{n}^{n}/c_{n}\rightarrow +\infty $. In the special case where $%
c_{n}=n^{\alpha }$, for $\alpha \in (0,1)$, $\lim\limits_{n\rightarrow
\infty }p_{n}/\rho _{n}^{2n}=0$ no matter whether $p_{n}=2$ or $\log n$ or $%
2\log \log n$. In this case, all the well-known information criteria can
consistently select the true model.
\end{remark}

\begin{theorem}
\label{thm02c}Let Assumption \ref{assIN} (i) holds. Assume the true DGP is
the EX model.

\begin{enumerate}
\item[(1)] When $\lim\limits_{n\rightarrow \infty }\dfrac{p_{n}}{\rho ^{2n}}%
=0,$ 
\begin{align*}
\lim\limits_{n\rightarrow \infty }P\left\{ \widehat{k}_{IC}=0|k=1\right\} &
=\lim\limits_{n\rightarrow \infty }P\left\{ \dfrac{n}{\rho ^{2n}}\left(
IC_{1}-IC_{0}\right) >0\right\} =0, \\
\lim\limits_{n\rightarrow \infty }P\left\{ \widehat{k}_{IC}=1|k=1\right\} &
=\lim\limits_{n\rightarrow \infty }P\left\{ \dfrac{n}{\rho ^{2n}}\left(
IC_{1}-IC_{0}\right) \leq 0\right\} =1.
\end{align*}

\item[(2)] When $\lim\limits_{n\rightarrow \infty }\dfrac{p_{n}}{\rho ^{2n}}%
=\pi \in (0,+\infty ),$ 
\begin{align*}
\lim\limits_{n\rightarrow \infty }P\left\{ \widehat{k}_{IC}=0|k=1\right\} &
=\lim\limits_{n\rightarrow \infty }P\left\{ \dfrac{n}{\rho ^{2n}}\left(
IC_{1}-IC_{0}\right) >0\right\} =P\left( \chi ^{2}(1)<(1+\rho )^{2}\pi
\right) , \\
\lim\limits_{n\rightarrow \infty }P\left\{ \widehat{k}_{IC}=1|k=1\right\} &
=\lim\limits_{n\rightarrow \infty }P\left\{ \dfrac{n}{\rho ^{2n}}\left(
IC_{1}-IC_{0}\right) \leq 0\right\} =1-P\left( \chi ^{2}(1)<(1+\rho )^{2}\pi
\right) .
\end{align*}

\item[(3)] When $\lim\limits_{n\rightarrow \infty }\dfrac{p_{n}}{\rho ^{2n}}%
\rightarrow +\infty ,$ 
\begin{align*}
\lim\limits_{n\rightarrow \infty }P\left\{ \widehat{k}_{IC}=0|k=1\right\} &
=\lim\limits_{n\rightarrow \infty }P\left\{ \dfrac{n}{p_{n}}\left(
IC_{1}-IC_{0}\right) >0\right\} =1, \\
\lim\limits_{n\rightarrow \infty }P\left\{ \widehat{k}_{IC}=1|k=1\right\} &
=\lim\limits_{n\rightarrow \infty }P\left\{ \dfrac{n}{p_{n}}\left(
IC_{1}-IC_{0}\right) \leq 0\right\} =0.
\end{align*}
\end{enumerate}
\end{theorem}

\begin{remark}
\label{rmk02c}Theorem \ref{thm02c} shows that the limit probability of
selecting the correct model by information criteria under the EX model
depends also critically on two parameters, $\rho $, $p_{n}$. As expected,
the larger $\rho $, the higher probability for the information criteria to
select the correct model. Interestingly, the smaller $p_{n}$, the higher
probability for the information criteria to select the correct model. If $%
p_{n}=2$ or $\log n$ or $2\log \log n$, $\lim\limits_{n\rightarrow \infty
}p_{n}/\rho ^{2n}=0$ and hence case (1) applies, suggesting that all the
well-known information criteria can consistently select the true model.
\end{remark}

Results in Theorem \ref{thm02b} can be extended to cover the LTUE model and
the ME model with weakly dependent errors. The following proposition
establishes the results for the ME model.

\begin{proposition}
\label{prop05}Let\textbf{\ }Assumption \ref{assIN} (i) or (ii) and the
assumption specified in Equation (\ref{dep}) hold. Assume the true DGP is
the ME model.

\begin{enumerate}
\item[(1)] When $\lim\limits_{n\rightarrow \infty }\dfrac{p_{n}}{\rho
_{n}^{2n}}=0,$ 
\begin{align*}
\lim\limits_{n\rightarrow \infty }P\left\{ \widehat{k}_{IC}=0|k=1\right\} &
=\lim\limits_{n\rightarrow \infty }P\left\{ \dfrac{n}{\rho _{n}^{2n}}\left(
IC_{1}-IC_{0}\right) >0\right\} =0, \\
\lim\limits_{n\rightarrow \infty }P\left\{ \widehat{k}_{IC}=1|k=1\right\} &
=\lim\limits_{n\rightarrow \infty }P\left\{ \dfrac{n}{\rho _{n}^{2n}}\left(
IC_{1}-IC_{0}\right) \leq 0\right\} =1.
\end{align*}

\item[(2)] When $\lim\limits_{n\rightarrow \infty }\dfrac{p_{n}}{\rho
_{n}^{2n}}=\pi \in (0,+\infty ),$ 
\begin{align*}
\lim\limits_{n\rightarrow \infty }P\left\{ \widehat{k}_{IC}=0|k=1\right\} &
=\lim\limits_{n\rightarrow \infty }P\left\{ \dfrac{n}{\rho _{n}^{2n}}\left(
IC_{1}-IC_{0}\right) >0\right\} =P\left( \chi ^{2}(1)<\frac{4\pi }{\omega
^{2}}\right) , \\
\lim\limits_{n\rightarrow \infty }P\left\{ \widehat{k}_{IC}=1|k=1\right\} &
=\lim\limits_{n\rightarrow \infty }P\left\{ \dfrac{n}{\rho _{n}^{2n}}\left(
IC_{1}-IC_{0}\right) \leq 0\right\} =1-P\left( \chi ^{2}(1)<\frac{4\pi }{%
\omega ^{2}}\right) .
\end{align*}

where $\omega ^{2}=\left( \sum\nolimits_{j=0}^{\infty }F_{j}\right) ^{2}.$

\item[(3)] When $\dfrac{p_{n}}{\rho _{n}^{2n}}\rightarrow +\infty ,$ 
\begin{align*}
\lim\limits_{n\rightarrow \infty }P\left\{ \widehat{k}_{IC}=0|k=1\right\} &
=\lim\limits_{n\rightarrow \infty }P\left\{ \dfrac{n}{p_{n}}\left(
IC_{1}-IC_{0}\right) >0\right\} =1, \\
\lim\limits_{n\rightarrow \infty }P\left\{ \widehat{k}_{IC}=1|k=1\right\} &
=\lim\limits_{n\rightarrow \infty }P\left\{ \dfrac{n}{p_{n}}\left(
IC_{1}-IC_{0}\right) \leq 0\right\} =0.
\end{align*}
\end{enumerate}
\end{proposition}

\section{Limit Properties Based on the Indirect Inference Estimator}

The OLS estimator of $\rho _{n}$ in Model (\ref{eqUR}) is known to be biased
and the bias is acute when $\rho _{n}$ is close to unity. To reduce the
bias, the indirect inference method of Smith (1993) and Gour\'{e}rioux et al
(1993) can be used if Model (\ref{eqUR}) is fully specified. Phillips (2012)
derives the asymptotic theory of the indirect inference estimator when the
model is UR or LTU and $u_{t}\overset{iid}{\sim }\mathcal{N}(0,\sigma ^{2})$%
. Throughout the paper we denote $\breve{\rho}$ the indirect inference
estimator of $\rho $. Let $h(c)=c+g(c)$ and $g(c)=g^{-}(c)1_{\{c\leq
0\}}+g^{+}(c)1_{\{c>0\}}$ with 
\begin{alignat*}{2}
g^{-}(c)& = & & -\dfrac{3}{4}\int_{0}^{\infty }e^{-\frac{v}{4}%
}k^{-}(v;c)^{1/2}dv+\dfrac{1}{4}\int_{0}^{\infty }e^{-\frac{v}{4}%
}k^{-}(v;c)^{3/2}dv \\
& {} & & -\dfrac{e^{2c}}{8}\int_{0}^{\infty }e^{-\frac{5v}{4}%
}k^{-}(v;c)^{3/2}vdv, \\
g^{+}(c)& = & & \dfrac{3}{4}\int_{0}^{\infty }e^{\frac{w}{4}%
}k^{+}(w;c)^{1/2}dw-\dfrac{1}{4}\int_{0}^{\infty }e^{\frac{w}{4}%
}k^{+}(w;c)^{3/2}dw \\
& {} & & -\dfrac{e^{2c}}{8}\int_{0}^{\infty }e^{\frac{5w}{4}%
}k^{+}(w;c)^{3/2}wdw, \\
k^{-}(v;c)& = & & \dfrac{2v-4c}{v+e^{2c}ve^{-v}-4c}, \\
k^{+}(w;c)& = & & \dfrac{2w+4c}{w+e^{2c}we^{w}+4c}.
\end{alignat*}%
Phillips (2012) shows that under the UR model,

\begin{equation*}
n\left( \breve{\rho}-1\right) \Rightarrow h^{-1}\left(
\int_{0}^{1}BdB/\int_{0}^{1}B^{2}\right) \text{ as }n\rightarrow +\infty ,
\end{equation*}%
and under the LTUE model,%
\begin{equation*}
n\left( \breve{\rho}-\rho _{n}\right) \Rightarrow h^{-1}\left(
\int_{0}^{1}J_{c}dB/\int_{0}^{1}J_{c}^{2}+c\right) -c\text{ as }n\rightarrow
+\infty .
\end{equation*}

Let $\breve{k}_{IC}=0$ or $1$ mean the information criterion of the UR\
model is smaller or larger than that of the competing model when the model
is estimated by the indirect inference method. We aim to find is the limit
of the following probabilities: 
\begin{align}
& \lim\limits_{n\rightarrow \infty }P\left\{ \breve{k}_{IC}=0|k=1\right\} ;
\\
& \lim\limits_{n\rightarrow \infty }P\left\{ \breve{k}_{IC}=1|k=1\right\} ;
\\
& \lim\limits_{n\rightarrow \infty }P\left\{ \breve{k}_{IC}=0|k=0\right\} ;
\\
& \lim\limits_{n\rightarrow \infty }P\left\{ \breve{k}_{IC}=1|k=0\right\} .
\end{align}

\begin{theorem}
\label{thm03}Under Assumption \ref{assIN}(i) or (ii) or (iii), we have

\begin{enumerate}
\item[(1)] when $p_{n}\rightarrow \infty $ and $p_{n}/n\rightarrow 0$ as $%
n\rightarrow \infty $, 
\begin{align*}
\lim\limits_{n\rightarrow \infty }P\left\{ \breve{k}_{IC}=0|k=0\right\} &
=\lim\limits_{n\rightarrow \infty }P\left\{ IC_{0}-IC_{1}\leq 0\right\} =1,
\\
\lim\limits_{n\rightarrow \infty }P\left\{ \breve{k}_{IC}=1|k=0\right\} &
=\lim\limits_{n\rightarrow \infty }P\left\{ IC_{0}-IC_{1}>0\right\} =0;
\end{align*}

\item[(2)] when $p_{n}=2$, the asymptotic distribution under the AIC
criterion is 
\begin{align*}
\lim\limits_{n\rightarrow \infty }P\left\{ \breve{k}_{AIC}=0|k=0\right\} &
=P\left( \varsigma ^{2}<2\right) , \\
\lim\limits_{n\rightarrow \infty }P\left\{ \breve{k}_{AIC}=1|k=0\right\} &
=1-P\left( \varsigma ^{2}<2\right) ,
\end{align*}%
where 
\begin{equation*}
\varsigma ^{2}=%
\begin{cases}
\int_{0}^{1}B^{2}\cdot h^{-1}\left( \left( \dfrac{\int_{0}^{1}BdB}{%
\int_{0}^{1}B^{2}}\right) ^{2}\right) -2\int_{0}^{1}BdB\cdot h^{-1}\left( 
\dfrac{\int_{0}^{1}BdB}{\int_{0}^{1}B^{2}}\right) , & \text{if }\tau =0 \\ 
\int_{0}^{1}B_{\tau }^{2}\cdot h^{-1}\left( \left( \dfrac{%
\int_{0}^{1}B_{\tau }dB}{\int_{0}^{1}B_{\tau }^{2}}\right) ^{2}\right)
-2\int_{0}^{1}B_{\tau }dB\cdot h^{-1}\left( \dfrac{\int_{0}^{1}B_{\tau }dB}{%
\int_{0}^{1}B_{\tau }^{2}}\right) , & \text{if }\tau \in (0,\infty ) \\ 
h^{-1}\left( \mathcal{C}\right) ^{2}B_{0}^{2}(1)-2h^{-1}\left( \mathcal{C}%
\right) B(1)B_{0}(1), & \text{if }\tau =\infty%
\end{cases}%
,
\end{equation*}%
with $\mathcal{C}$ being a standard Cauchy variate.
\end{enumerate}
\end{theorem}

\begin{remark}
\label{rmk02}According to Theorem \ref{thm03}, as long as $p_{n}\rightarrow
\infty $ and $p_{n}/n\rightarrow 0$, information criteria based on the
indirect inference estimator is consistent in selecting the UR\ model.
Hence, BIC and HQIC based on the indirect inference estimator can
consistently select the UR model. Like the AIC criterion that is based on
the OLS estimator, the AIC criterion based on the indirect inference
estimator continues to be inconsistent. However, its asymptotic distribution
depends on $\varsigma ^{2}$, the squared unit root $t$-statistic for the
indirect inference estimator.
\end{remark}

\begin{remark}
\label{rmk03}As shown in Phillips (2012), the squared unit root $t$%
-statistic for the indirect inference estimator has a smaller variance than
that of the squared unit root $t$-statistic for the OLS estimator.
Consequently, $P\left( \varsigma ^{2}<2\right) >P\left( \xi ^{2}<2\right) $,
suggesting that AIC based on the indirect inference estimator can select the
true model (i.e. the UR model) with a larger probability than that based on
the OLS estimator.
\end{remark}

\begin{theorem}
\label{thm04a}Let Assumption \ref{assIN} (i) or (ii) holds. Assume the true
DGP is the LTUE model.

\begin{enumerate}
\item[(1)] When $p_{n}\rightarrow \infty $ and $p_{n}/n\rightarrow 0$ as $%
n\rightarrow \infty $, 
\begin{align*}
\lim\limits_{n\rightarrow \infty }P\left\{ \breve{k}_{IC}=0|k=1\right\} &
=\lim\limits_{n\rightarrow \infty }P\left\{ \dfrac{n}{p_{n}}\left(
IC_{1}-IC_{0}\right) >0\right\} =1, \\
\lim\limits_{n\rightarrow \infty }P\left\{ \breve{k}_{IC}=1|k=1\right\} &
=\lim\limits_{n\rightarrow \infty }P\left\{ \dfrac{n}{p_{n}}\left(
IC_{1}-IC_{0}\right) \leq 0\right\} =0.
\end{align*}

\item[(2)] When $p_{n}=2$, the asymptotic distribution under the AIC
criterion is 
\begin{align*}
\lim\limits_{n\rightarrow \infty }P\left\{ \breve{k}_{\text{AIC}%
}=0|k=1\right\} & =\lim\limits_{n\rightarrow \infty }P\left\{ {n}\left(
AIC_{1}-AIC_{0}\right) >0\right\} =1-P\left( \vartheta ^{2}>2\right) , \\
\lim\limits_{n\rightarrow \infty }P\left\{ \breve{k}_{\text{AIC}%
}=1|k=1\right\} & =\lim\limits_{n\rightarrow \infty }P\left\{ {n}\left(
AIC_{1}-AIC_{0}\right) \leq 0\right\} =P\left( \vartheta ^{2}>2\right) ,
\end{align*}%
where 
\begin{equation*}
\vartheta ^{2}\equiv 2h^{-1}\left( \dfrac{\int_{0}^{1}J_{c}dB}{%
\int_{0}^{1}J_{c}^{2}}+c\right) \left(
\int_{0}^{1}J_{c}dB+c\int_{0}^{1}J_{c}^{2}\right) -h^{-1}\left( \dfrac{%
\int_{0}^{1}J_{c}dB}{\int_{0}^{1}J_{c}^{2}}+c\right)
^{2}\int_{0}^{1}J_{c}^{2}.
\end{equation*}
\end{enumerate}
\end{theorem}

\begin{remark}
\label{rmk04a}Theorem \ref{thm04a} shows that all the information criteria
continue to be inconsistent in distinguishing between the LTUE model and the
UR models when data come from the LTUE model even when the indirect
inference estimation is employed. AIC selects the wrong model with
probability going to $1-P\left( \vartheta ^{2}>2\right) $. Since the
variance of $\zeta ^{2}$ is bigger than that of $\vartheta ^{2}$, the tail
probability of $\zeta ^{2}$ is larger than that of $\vartheta ^{2}$,
suggesting that AIC based on OLS selects the true model (i.e. LTUE model)
with a greater\ probability than AIC\ based on the indirect inference
estimator. This is a rather surprising result and suggests that the
superiority in estimation does not necessarily translate to the superiority
in model selection.
\end{remark}

\begin{theorem}
\label{thm04b}Let Assumption \ref{assIN} (i) or (ii) holds. Assume the true
DGP is the ME model.

\begin{enumerate}
\item[(1)] When $\lim\limits_{n\rightarrow \infty }\dfrac{p_{n}}{\rho
_{n}^{2n}}=0,$ 
\begin{align*}
\lim\limits_{n\rightarrow \infty }P\left\{ \breve{k}_{IC}=0|k=1\right\} &
=\lim\limits_{n\rightarrow \infty }P\left\{ \dfrac{n}{\rho _{n}^{2n}}\left(
IC_{1}-IC_{0}\right) >0\right\} =0, \\
\lim\limits_{n\rightarrow \infty }P\left\{ \breve{k}_{IC}=1|k=1\right\} &
=\lim\limits_{n\rightarrow \infty }P\left\{ \dfrac{n}{\rho _{n}^{2n}}\left(
IC_{1}-IC_{0}\right) \leq 0\right\} =1.
\end{align*}

\item[(2)] When $\lim\limits_{n\rightarrow \infty }\dfrac{p_{n}}{\rho
_{n}^{2n}}=\pi \in (0,+\infty ),$ 
\begin{align*}
\lim\limits_{n\rightarrow \infty }P\left\{ \breve{k}_{IC}=0|k=1\right\} &
=\lim\limits_{n\rightarrow \infty }P\left\{ \dfrac{n}{\rho _{n}^{2n}}\left(
IC_{1}-IC_{0}\right) >0\right\} =P\left( \chi ^{2}(1)<4\pi \right) , \\
\lim\limits_{n\rightarrow \infty }P\left\{ \breve{k}_{IC}=1|k=1\right\} &
=\lim\limits_{n\rightarrow \infty }P\left\{ \dfrac{n}{\rho _{n}^{2n}}\left(
IC_{1}-IC_{0}\right) \leq 0\right\} =1-P\left( \chi ^{2}(1)<4\pi \right) .
\end{align*}

\item[(3)] When $\dfrac{p_{n}}{\rho _{n}^{2n}}\rightarrow +\infty ,$ 
\begin{align*}
\lim\limits_{n\rightarrow \infty }P\left\{ \breve{k}_{IC}=0|k=1\right\} &
=\lim\limits_{n\rightarrow \infty }P\left\{ \dfrac{n}{p_{n}}\left(
IC_{1}-IC_{0}\right) >0\right\} =1, \\
\lim\limits_{n\rightarrow \infty }P\left\{ \breve{k}_{IC}=1|k=1\right\} &
=\lim\limits_{n\rightarrow \infty }P\left\{ \dfrac{n}{p_{n}}\left(
IC_{1}-IC_{0}\right) \leq 0\right\} =0.
\end{align*}
\end{enumerate}
\end{theorem}

\begin{remark}
\label{rmk04b}The results in Theorem \ref{thm04b} are the same as those in
Theorem \ref{thm02b}, suggesting all the well-known information criteria can
consistently select the true model (i.e. ME model) when $c_{n}=n^{\alpha }$,
for $\alpha \in (0,1)$.
\end{remark}

\begin{theorem}
\label{thm04c}Let Assumption \ref{assIN} (i) holds. Assume the true DGP is
the EX model.

\begin{enumerate}
\item[(1)] When $\lim\limits_{n\rightarrow \infty }\dfrac{p_{n}}{\rho ^{2n}}%
=0,$ 
\begin{align*}
\lim\limits_{n\rightarrow \infty }P\left\{ \breve{k}_{IC}=0|k=1\right\} &
=\lim\limits_{n\rightarrow \infty }P\left\{ \dfrac{n}{\rho ^{2n}}\left(
IC_{1}-IC_{0}\right) >0\right\} =0, \\
\lim\limits_{n\rightarrow \infty }P\left\{ \breve{k}_{IC}=1|k=1\right\} &
=\lim\limits_{n\rightarrow \infty }P\left\{ \dfrac{n}{\rho ^{2n}}\left(
IC_{1}-IC_{0}\right) \leq 0\right\} =1.
\end{align*}

\item[(2)] When $\lim\limits_{n\rightarrow \infty }\dfrac{p_{n}}{\rho ^{2n}}%
=\pi \in (0,+\infty ),$ 
\begin{align*}
\lim\limits_{n\rightarrow \infty }P\left\{ \breve{k}_{IC}=0|k=1\right\} &
=\lim\limits_{n\rightarrow \infty }P\left\{ \dfrac{n}{\rho ^{2n}}\left(
IC_{1}-IC_{0}\right) >0\right\} =P\left( \chi ^{2}(1)<(1+\rho )^{2}\pi
\right) , \\
\lim\limits_{n\rightarrow \infty }P\left\{ \breve{k}_{IC}=1|k=1\right\} &
=\lim\limits_{n\rightarrow \infty }P\left\{ \dfrac{n}{\rho ^{2n}}\left(
IC_{1}-IC_{0}\right) \leq 0\right\} =1-P\left( \chi ^{2}(1)<(1+\rho )^{2}\pi
\right) .
\end{align*}

\item[(3)] When $\lim\limits_{n\rightarrow \infty }\dfrac{p_{n}}{\rho ^{2n}}%
\rightarrow +\infty ,$ 
\begin{align*}
\lim\limits_{n\rightarrow \infty }P\left\{ \breve{k}_{IC}=0|k=1\right\} &
=\lim\limits_{n\rightarrow \infty }P\left\{ \dfrac{n}{p_{n}}\left(
IC_{1}-IC_{0}\right) >0\right\} =1, \\
\lim\limits_{n\rightarrow \infty }P\left\{ \breve{k}_{IC}=1|k=1\right\} &
=\lim\limits_{n\rightarrow \infty }P\left\{ \dfrac{n}{p_{n}}\left(
IC_{1}-IC_{0}\right) \leq 0\right\} =0.
\end{align*}
\end{enumerate}
\end{theorem}

\begin{remark}
\label{rmk04c} The results in Theorem \ref{thm04c} are the same as those in
Theorem \ref{thm02c}, suggesting that all the well-known information
criteria can consistently select the true model (i.e. EX model).
\end{remark}

\section{Monte Carlo Study}

\label{secSim}

In this section, we examine the performance of alternative information
criteria, namely, AIC, BIC and HQIC, in finite sample via simulated data and
check the reliability of the asymptotic results developed in Section 3 and
Section 4. In the simulation study, we use both OLS and the indirect
inference method to estimate $\rho _{n}$ from sample paths that are
simulated from different DGPs. In total we design four experiments. In the
first experiment we simulate data from the UR model. In the second
experiment we simulate data from the LTUE model with $c=1$ (i.e. $\rho
_{n}=1+1/n)$. In the third experiment we simulate data from two ME models
with $c_{n}=n^{0.1}$, $n^{0.3}$, respectively. In the last experiment we
simulate data from the EX model with $\rho =1.01,1.05$, respectively. In all
experiments, we simulate 10,000 sample paths with initial value $X_{0}=0$
and four sample sizes are considered, $n=100,200,500,1000$. In each
experiment, we report the fraction of the number of times in which the
correct model is selected out of 10,000 replications.

Table \ref{tab01} reports the results when the true DGP is UR. Several
results can be found here. First, the probability for BIC\ and HQIC to
select the true model grows as $n$ grows. However, the probability for AIC
to select the true model does not seem to increase or decrease as $n$ grows.
This observation is consistent with the asymptotic results reported in
Theorem \ref{thm01}. Second, the probability for BIC to select the true
model is larger than that in HQIC which is in turn larger than AIC in these
four sample sizes. So we can conclude that the probability grows as $p_{n}$
increases since $2<2\log \log n<\log n$ when $100\leq n\leq 1000$. Third,
the probability implied by AIC based on the indirect inference estimator is
larger than that based on OLS. This finding is consistent with Theorem \ref%
{thm03} and Remark \ref{rmk03}.

\begin{table}[tph]
\caption{Probability of Selecting the Correct Model when Data Come from the
UR Model}
\label{tab01}
\begin{center}
\begin{tabular}{l|ccc|ccc}
\hline\hline
$n$ & \multicolumn{3}{c|}{100} & \multicolumn{3}{c}{200} \\ \hline
IC & AIC & BIC & HQIC & AIC & BIC & HQIC \\ 
OLS & 0.8160 & 0.9604 & 0.9020 & 0.8155 & 0.9751 & 0.9249 \\ 
IIE & 0.8731 & 0.9702 & 0.9292 & 0.8742 & 0.9810 & 0.9445 \\ \hline
$n$ & \multicolumn{3}{c|}{500} & \multicolumn{3}{c}{1000} \\ \hline
IC & AIC & BIC & HQIC & AIC & BIC & HQIC \\ 
OLS & 0.8127 & 0.9849 & 0.9335 & 0.8195 & 0.9895 & 0.9402 \\ 
IIE & 0.8704 & 0.9881 & 0.9508 & 0.8759 & 0.9918 & 0.9566 \\ \hline\hline
\end{tabular}%
\end{center}
\end{table}

Table \ref{tab02} report the results when the true DGP is the LTUE model
with $c_{n}=1$. Also reported is the value of $p_{n}/\rho _{n}^{2n}$.
Several results can be found here. First, the probability for BIC\ and HQIC
to select the true model becomes smaller as $n$ grows. However, the
probability for AIC to select the true model does not seem to increase or
decrease as $n$ grows. This observation is consistent with the asymptotic
results in Theorem \ref{thm02a}. Second, the probability implied by AIC
based on the indirect inference estimator is smaller than that based on OLS.
This finding is consistent with in Theorem \ref{thm04a} and Remark \ref%
{rmk04a}. Finally, it seems that AIC performs better than BIC and HQIC in
all cases.

\begin{table}[tph]
\caption{Probability of Selecting the Correct Model when Data Come from the
LTUE Model with $c_{n}=1$}
\label{tab02}
\begin{center}
\begin{tabular}{c|ccc|ccc}
\hline\hline
$n$ & \multicolumn{3}{c|}{100} & \multicolumn{3}{c}{200} \\ \hline
\multicolumn{1}{l|}{IC} & AIC & BIC & HQIC & AIC & BIC & HQIC \\ 
\multicolumn{1}{l|}{$p_{n}/\rho _{n}^{2n}$} & 0.2734 & 0.6295 & 0.4175 & 
0.2720 & 0.7206 & 0.4536 \\ 
\multicolumn{1}{l|}{OLS} & 0.3516 & 0.1475 & 0.2420 & 0.3406 & 0.1305 & 
0.2156 \\ 
\multicolumn{1}{l|}{IIE} & 0.1485 & 0.0445 & 0.0922 & 0.1235 & 0.0269 & 
0.0663 \\ \hline
$n$ & \multicolumn{3}{c|}{500} & \multicolumn{3}{c}{1000} \\ \hline
\multicolumn{1}{l|}{IC} & AIC & BIC & HQIC & AIC & BIC & HQIC \\ 
\multicolumn{1}{l|}{$p_{n}/\rho _{n}^{2n}$} & 0.2712 & 0.8427 & 0.4955 & 
0.2709 & 0.9358 & 0.5236 \\ 
\multicolumn{1}{l|}{OLS} & 0.3474 & 0.1019 & 0.1933 & 0.3416 & 0.0871 & 
0.1823 \\ 
\multicolumn{1}{l|}{IIE} & 0.1169 & 0.0134 & 0.0517 & 0.1089 & 0.0090 & 
0.0394 \\ \hline\hline
\end{tabular}%
\end{center}
\end{table}

Table \ref{tab02b} report the results when the true DGP is the ME model with 
$c_{n}=n^{0.1},n^{0.3}$. Also reported is the value of $p_{n}/\rho _{n}^{2n}$%
. Several results can be found here. First, the probability for all three
information criteria to select the true model grows as $n$ increases. This
observation is consistent with the asymptotic results reported in Theorem %
\ref{thm02b} and Remark \ref{rmk04b}. Second, comparing the results for $%
c_{n}=n^{0.1}$ and those for $c_{n}=n^{0.3}$, the probability for all three
information criteria to select the true model increases when $c_{n}$ is
bigger. Third, the probability based on the indirect inference estimator is
smaller than that based on OLS. Finally, it seems that AIC performs better
than BIC and HQIC in all cases.

\begin{table}[tph]
\caption{Probability of Selecting the Correct Model when Data Come from the
ME Model with $c_{n}=n^{0.1}$ and $c_{n}=n^{0.3}$}
\label{tab02b}
\begin{center}
\begin{tabular}{c|ccc|ccc}
\hline\hline
\multicolumn{7}{c}{ME Model with $c_{n}=n^{0.1}$} \\ \hline
$n$ & \multicolumn{3}{c|}{100} & \multicolumn{3}{c}{200} \\ \hline
\multicolumn{1}{l|}{IC} & AIC & BIC & HQIC & AIC & BIC & HQIC \\ 
\multicolumn{1}{l|}{$p_{n}/\rho _{n}^{2n}$} & 0.0861 & 0.1983 & 0.1316 & 
0.0679 & 0.1799 & 0.1132 \\ 
\multicolumn{1}{l|}{OLS} & 0.5183 & 0.3403 & 0.4349 & 0.5554 & 0.3638 & 
0.4629 \\ 
\multicolumn{1}{l|}{IIE} & 0.3071 & 0.1741 & 0.2406 & 0.3211 & 0.1624 & 
0.2250 \\ \hline
$n$ & \multicolumn{3}{c|}{500} & \multicolumn{3}{c}{1000} \\ \hline
\multicolumn{1}{l|}{IC} & AIC & BIC & HQIC & AIC & BIC & HQIC \\ 
\multicolumn{1}{l|}{$p_{n}/\rho _{n}^{2n}$} & 0.0486 & 0.1512 & 0.0889 & 
0.0371 & 0.1282 & 0.0718 \\ 
\multicolumn{1}{l|}{OLS} & 0.6151 & 0.4083 & 0.5048 & 0.6469 & 0.4374 & 
0.5494 \\ 
\multicolumn{1}{l|}{IIE} & 0.3544 & 0.2008 & 0.2815 & 0.3925 & 0.2351 & 
0.3129 \\ \hline\hline
\multicolumn{7}{c}{ME Model with $c_{n}=n^{0.3}$} \\ \hline
$n$ & \multicolumn{3}{c|}{100} & \multicolumn{3}{c}{200} \\ \hline
\multicolumn{1}{l|}{IC} & AIC & BIC & HQIC & AIC & BIC & HQIC \\ 
\multicolumn{1}{l|}{$p_{n}/\rho _{n}^{2n}$} & 0.0008 & 0.0019 & 0.0012 & 
0.0001 & 0.0003 & 0.0002 \\ 
\multicolumn{1}{l|}{OLS} & 0.9374 & 0.9066 & 0.9235 & 0.9749 & 0.9608 & 
0.9683 \\ 
\multicolumn{1}{l|}{IIE} & 0.9274 & 0.8979 & 0.9163 & 0.9716 & 0.9578 & 
0.9648 \\ \hline
$n$ & \multicolumn{3}{c|}{500} & \multicolumn{3}{c}{1000} \\ \hline
\multicolumn{1}{l|}{IC} & AIC & BIC & HQIC & AIC & BIC & HQIC \\ 
\multicolumn{1}{l|}{$p_{n}/\rho _{n}^{2n}$} & 1.0e-06 & 1.0e-05 & 1.0e-06 & 
1.0e-07 & 1.0e-07 & 1.0e-07 \\ 
\multicolumn{1}{l|}{OLS} & 0.9948 & 0.9907 & 0.9938 & 0.9988 & 0.9985 & 
0.9986 \\ 
\multicolumn{1}{l|}{IIE} & 0.9938 & 0.9901 & 0.9933 & 0.9986 & 0.9985 & 
0.9985 \\ \hline\hline
\end{tabular}%
\end{center}
\end{table}

Table \ref{tab03} report the results when the true DGP is the EX model with $%
\rho =1.01,1.05$. Also reported is the value of $p_{n}/\rho ^{2n}$. Several
results can be found here. First, when $\rho =1.01$, which is larger than
the unity by 1\%, the probability for information criteria to select the
correct model is small in all cases when the sample size is small. However,
it grows very quickly with the sample size. When $\rho =1.05$, the
probability for information criteria to select the correct model is almost 1
in all cases even when the sample size is small and increases with the
sample size. Finally, it seems that AIC performs better than BIC and HQIC in
all cases. 
\begin{table}[tph]
\caption{Probability of Selecting the Correct Model when Data Come from the
Regular Explosive Model with $\protect\rho =1.01,1.05$.}
\label{tab03}
\begin{center}
\begin{tabular}{c|ccc|ccc}
\hline\hline
\multicolumn{7}{c}{Explosive Model with $\rho =1.01$} \\ \hline
$n$ & \multicolumn{3}{c|}{100} & \multicolumn{3}{c}{200} \\ \hline
\multicolumn{1}{l|}{IC} & AIC & BIC & HQIC & AIC & BIC & HQIC \\ 
\multicolumn{1}{l|}{$p_{n}/\rho ^{2n}$} & 0.2734 & 0.6295 & 0.4175 & 0.0374
& 0.0990 & 0.0623 \\ 
\multicolumn{1}{l|}{OLS} & 0.3516 & 0.1475 & \multicolumn{1}{c|}{0.2420} & 
0.6449 & 0.4820 & 0.5555 \\ 
\multicolumn{1}{l|}{IIE} & 0.1485 & 0.0445 & \multicolumn{1}{c|}{0.0922} & 
0.4740 & 0.3059 & 0.3845 \\ \hline
$n$ & \multicolumn{3}{c|}{500} & \multicolumn{3}{c}{1000} \\ \hline
\multicolumn{1}{l|}{IC} & AIC & BIC & HQIC & AIC & BIC & HQIC \\ 
\multicolumn{1}{l|}{$p_{n}/\rho ^{2n}$} & 1.0e-4 & 1.0e-4 & 1.0e-4 & 1.0e-9
& 1.0e-8 & 1.0e-9 \\ 
\multicolumn{1}{l|}{OLS} & 0.9775 & 0.9599 & 0.9704 & 0.9998 & 0.9997 & 
0.9998 \\ 
\multicolumn{1}{l|}{IIE} & 0.9733 & 0.9563 & 0.9681 & 0.9998 & 0.9997 & 
0.9998 \\ \hline\hline
\multicolumn{7}{c}{Explosive Model with $\rho =1.05$} \\ \hline
$n$ & \multicolumn{3}{c|}{100} & \multicolumn{3}{c}{200} \\ \hline
\multicolumn{1}{l|}{IC} & AIC & BIC & HQIC & AIC & BIC & HQIC \\ 
\multicolumn{1}{l|}{$p_{n}/\rho ^{2n}$} & 0.0001 & 0.0003 & 0.0002 & 1.0e-07
& 1.0e-07 & 1.0e-07 \\ 
\multicolumn{1}{l|}{OLS} & 0.9741 & 0.9643 & 0.9681 & 0.9999 & 0.9998 & 
0.9998 \\ 
\multicolumn{1}{l|}{IIE} & 0.9703 & 0.9626 & 0.9655 & 0.9999 & 0.9998 & 
0.9998 \\ \hline
$n$ & \multicolumn{3}{c|}{500} & \multicolumn{3}{c}{1000} \\ \hline
\multicolumn{1}{l|}{IC} & AIC & BIC & HQIC & AIC & BIC & HQIC \\ 
\multicolumn{1}{l|}{$p_{n}/\rho ^{2n}$} & 1.0e-20 & 1.0e-20 & 1.0e-20 & 
1.0e-41 & 1.0e-41 & 1.0e-41 \\ 
\multicolumn{1}{l|}{OLS} & 1.0000 & 1.0000 & 1.0000 & 1.0000 & 1.0000 & 
1.0000 \\ 
\multicolumn{1}{l|}{IIE} & 1.0000 & 1.0000 & 1.0000 & 1.0000 & 1.0000 & 
1.0000 \\ \hline\hline
\end{tabular}%
\end{center}
\end{table}

\section{Conclusion}

\label{secCon}

This paper studies the limit properties of information criteria for
distinguishing between unit root model and three types of explosive models.
Both the OLS estimator and the indirect inference estimator are employed to
estimate the AR coefficient in the candidate model. This paper contributes
to the literature in three aspects. First, our results extends results in
the literature to the explosive side of the unit root, and we find that
information criteria consistently choose the unit root model when the unit
root model is the true model. Second, we show that the limiting
probabilities for information criteria to select the explosive model depends
on both the distance of autoregressive coefficient from unity and the size
of penalty term in the information criteria. When the penalty term is not
too large and the root is not too close to unit root, all the information
criteria consistently select the true model. It is known that the indirect
inference method is effective in reducing the bias in OLS estimation in all
cases as well as reducing the variance in OLS estimation in the UR\ model
and in the LTU model. However, when information criteria are used in
connection with the indirect inference estimation, the limiting
probabilities for information criteria to select the correct model can go up
or down relative to that with the OLS estimation, depending on the true DGP.
When the true DGP is the UR\ model, the indirect inference estimation
increases the probability. When the true DGP is the LTUE\ model or the ME
model or the EX model, the indirect inference estimation decreases the
probability. This rather surprising result suggests that the superiority in
estimation does not necessarily translate to the superiority in model
selection.

\section*{Appendix}

\label{secApp} \addcontentsline{toc}{section}{Appendices} %
\renewcommand{\thesubsection}{\Alph{subsection}}

\subsection{Proof of Theorem \protect\ref{thm01}}

The proof is same as the proof for Theorem 1 in Phillips (2008), and hence
omitted.

\subsection{Proof of Theorem \protect\ref{thm02a}}

When the true DGP is the LTUE model, we have $0<c<\infty $ and%
\begin{align*}
IC_{0}& =\log \widehat{\sigma }_{0}^{2}=\log \left\{ \frac{1}{n}%
\sum_{t=1}^{n}\left( X_{t}-X_{t-1}\right) ^{2}\right\} \\
& =\log \left\{ \frac{1}{n}\sum_{t=1}^{n}\left[ \left( \rho _{n}-1\right)
X_{t-1}+u_{t}\right] ^{2}\right\} \\
& =\log \left\{ \frac{1}{n}\left( \rho _{n}-1\right)
^{2}\sum_{t=1}^{n}X_{t-1}^{2}+\frac{2}{n}\left( \rho _{n}-1\right)
\sum_{t=1}^{n}X_{t-1}u_{t}+\frac{1}{n}\sum_{t=1}^{n}u_{t}^{2}\right\} .
\end{align*}

By Lemma 1 in Phillips (1987b), when the process is initialized at $X_{0}$,
we know 
\begin{equation}
\frac{1}{n^{2}}\sum_{t=1}^{n}X_{t-1}^{2}\Rightarrow \sigma
^{2}\int_{0}^{1}J_{c}^{2},  \label{eqLU01}
\end{equation}%
and 
\begin{equation}
\frac{1}{n}\sum_{t=1}^{n}X_{t-1}u_{t}\Rightarrow \sigma
^{2}\int_{0}^{1}J_{c}dB,  \label{eqLU02}
\end{equation}%
where%
\begin{equation*}
J_{c}(r)=\int_{0}^{r}e^{c(r-s)}dB(s).
\end{equation*}

Therefore, by Equation (\ref{eqLU01}) and (\ref{eqLU02}) we have 
\begin{align}
IC_{0}& =\log \left\{ \frac{\sigma ^{2}c^{2}}{n}\int_{0}^{1}J_{c}^{2}+\frac{%
2c\sigma ^{2}}{n}\int_{0}^{1}J_{c}dB+\sigma ^{2}+o_{p}(n^{-1})\right\} 
\notag \\
& =\log \sigma ^{2}+\log \left\{ 1+\frac{2c}{n}\int_{0}^{1}J_{c}(r)dB+\frac{%
c^{2}}{n}\int_{0}^{1}J_{c}^{2}+o_{p}(n^{-1})\right\} .  \label{eq01}
\end{align}

We also know from Phillips (1987b) that 
\begin{equation}
n\left( \widehat{\rho }_{n}-\rho _{n}\right) \Rightarrow \dfrac{%
\int_{0}^{1}J_{c}dB}{\int_{0}^{1}J_{c}^{2}},  \label{eqLU03}
\end{equation}%
Hence, 
\begin{align}
IC_{1}& =\log \widehat{\sigma }_{1}^{2}+\frac{p_{n}}{n}  \notag \\
& =\log \left\{ n^{-1}\sum_{t=1}^{n}\left( X_{t}-\widehat{\rho }%
_{n}X_{t-1}\right) ^{2}\right\} +\frac{p_{n}}{n}  \notag \\
& =\log \left\{ \frac{1}{n}\sum_{t=1}^{n}\left[ \left( \rho _{n}-\widehat{%
\rho }_{n}\right) X_{t-1}+u_{t}\right] ^{2}\right\} +\frac{p_{n}}{n}  \notag
\\
& =\log \left\{ \frac{1}{n}\left( \rho _{n}-\widehat{\rho }_{n}\right)
^{2}\sum_{t=1}^{n}X_{t-1}^{2}+\frac{2}{n}\left( \rho _{n}-\widehat{\rho }%
_{n}\right) \sum_{t=1}^{n}X_{t-1}u_{t}+\frac{1}{n}\sum_{t=1}^{n}u_{t}^{2}%
\right\} +\frac{p_{n}}{n}  \notag \\
& \Rightarrow \log \left\{ -\dfrac{\sigma ^{2}}{n}\dfrac{\left(
\int_{0}^{1}J_{c}dB\right) ^{2}}{\int_{0}^{1}J_{c}^{2}}+\sigma ^{2}\right\} +%
\frac{p_{n}}{n}  \notag \\
& =\log \sigma ^{2}+\log \left\{ 1-\dfrac{1}{n}\dfrac{\left(
\int_{0}^{1}J_{c}dB\right) ^{2}}{\int_{0}^{1}J_{c}^{2}}\right\} +\frac{p_{n}%
}{n}.  \label{eq02}
\end{align}

Therefore, by Equation (\ref{eq01}) and (\ref{eq02}), we have 
\begin{equation*}
IC_{1}-IC_{0}\Rightarrow \log \left\{ 1-\dfrac{1}{n}\dfrac{\left(
\int_{0}^{1}J_{c}dB\right) ^{2}}{\int_{0}^{1}J_{c}^{2}}\right\} -\log
\left\{ 1+\frac{2c}{n}\int_{0}^{1}J_{c}dB+\frac{c^{2}}{n}%
\int_{0}^{1}J_{c}^{2}\right\} +\frac{p_{n}}{n}.
\end{equation*}%
Hence, if $p_{n}=2$ (as in AIC), as $n\rightarrow \infty $, we have 
\begin{equation}
n\left( IC_{1}-IC_{0}\right) \Rightarrow 2-\dfrac{\left(
\int_{0}^{1}J_{c}dB\right) ^{2}}{\int_{0}^{1}J_{c}^{2}}-2c%
\int_{0}^{1}J_{c}dB-{c}^{2}\int_{0}^{1}J_{c}^{2}.
\end{equation}%
If $p_{n}\rightarrow \infty $ and $\dfrac{p_{n}}{n}\rightarrow 0$, we have 
\begin{equation*}
\dfrac{n}{p_{n}}\left( IC_{1}-IC_{0}\right) \Rightarrow 1.
\end{equation*}

\subsection{Proof of Theorem \protect\ref{thm02b}}

When the true DGP is the ME model, we have 
\begin{align}
IC_{0}& =\log \widehat{\sigma }_{0}^{2}=\log \left\{ \frac{1}{n}%
\sum_{t=1}^{n}\left( X_{t}-X_{t-1}\right) ^{2}\right\}  \notag \\
& =\log \left\{ \frac{1}{n}\sum_{t=1}^{n}\left[ \left( \rho _{n}-1\right)
X_{t-1}+u_{t}\right] ^{2}\right\}  \notag \\
& =\log \left\{ \frac{1}{n}\left( \rho _{n}-1\right)
^{2}\sum_{t=1}^{n}X_{t-1}^{2}+\frac{2}{n}\left( \rho _{n}-1\right)
\sum_{t=1}^{n}X_{t-1}u_{t}+\frac{1}{n}\sum_{t=1}^{n}u_{t}^{2}\right\} .
\end{align}

According to Phillips and Magdalinos (2007), when the process is initialized
at $X_{0}=o_{p}(\sqrt{n/c_{n}})$, we have 
\begin{equation}
\frac{c_{n}^{2}\rho _{n}^{-2n}}{n^{2}}\sum_{t=1}^{n}X_{t-1}^{2}\Rightarrow 
\frac{\sigma ^{2}}{4}Y^{2},  \label{eqPM01}
\end{equation}%
\begin{equation}
\frac{c_{n}\rho _{n}^{-n}}{n}\sum_{t=1}^{n}X_{t-1}u_{t}\Rightarrow \frac{%
\sigma ^{2}}{2}XY,  \label{eqPM02}
\end{equation}%
and 
\begin{equation}
\frac{n\rho _{n}^{n}}{2c_{n}}\left( \widehat{\rho }_{n}-\rho _{n}\right)
\Rightarrow \frac{X}{Y}\sim \mathcal{C},  \label{eqPM03}
\end{equation}%
where $X,Y\sim \mathcal{N}\left( 0,1\right) $ and $\mathcal{C}$ is a
standard Cauchy variate.

Therefore, by (\ref{eqPM01}) and (\ref{eqPM02}) we have 
\begin{align}
IC_{0}& \Rightarrow \log \left\{ \frac{\sigma ^{2}}{4n\rho _{n}^{-2n}}Y^{2}+%
\frac{\sigma ^{2}}{n\rho _{n}^{-n}}XY+\sigma ^{2}\right\}  \notag \\
& =\log \sigma ^{2}+\log \left\{ \frac{1}{4n\rho _{n}^{-2n}}Y^{2}+\frac{1}{%
n\rho _{n}^{-n}}XY+1\right\} .  \label{eq03}
\end{align}

On the other hand, 
\begin{align*}
IC_{1}& =\log \widehat{\sigma }_{1}^{2}+\frac{p_{n}}{n} \\
& =\log \left\{ n^{-1}\sum_{t=1}^{n}\left( X_{t}-\widehat{\rho }%
_{n}X_{t-1}\right) ^{2}\right\} +\frac{p_{n}}{n} \\
& =\log \left\{ \frac{1}{n}\sum_{t=1}^{n}\left[ \left( \rho _{n}-\widehat{%
\rho }_{n}\right) X_{t-1}+u_{t}\right] ^{2}\right\} +\frac{p_{n}}{n} \\
& =\log \left\{ \frac{1}{n}\left( \rho _{n}-\widehat{\rho }_{n}\right)
^{2}\sum_{t=1}^{n}X_{t-1}^{2}+\frac{2}{n}\left( \rho _{n}-\widehat{\rho }%
_{n}\right) \sum_{t=1}^{n}X_{t-1}u_{t}+\frac{1}{n}\sum_{t=1}^{n}u_{t}^{2}%
\right\} +\frac{p_{n}}{n}.
\end{align*}

By equation (\ref{eqPM01}) to (\ref{eqPM03}), we obtain 
\begin{align}
IC_{1}& \Rightarrow \log \left\{ \frac{1}{n}\frac{4c_{n}^{2}}{n^{2}\rho
_{n}^{2n}}\mathcal{C}^{2}\frac{n^{2}\sigma ^{2}}{4c_{n}^{2}\rho _{n}^{-2n}}%
Y^{2}-\frac{2}{n}\frac{2c_{n}}{n\rho _{n}^{n}}\mathcal{C}\frac{n\sigma ^{2}}{%
2c_{n}\rho _{n}^{-n}}XY+\sigma ^{2}\right\} +\frac{p_{n}}{n}  \notag \\
& =\log \sigma ^{2}+\log \left\{ -\frac{1}{n}X^{2}+1\right\} +\frac{p_{n}}{n}%
.  \label{eq04}
\end{align}%
Therefore, by equation (\ref{eq03}) and (\ref{eq04}), we have 
\begin{equation*}
IC_{1}-IC_{0}\Rightarrow \log \left\{ 1-\frac{1}{n}X^{2}\right\} -\log
\left\{ 1+\frac{1}{4n\rho _{n}^{-2n}}Y^{2}+\frac{1}{n\rho _{n}^{-n}}%
XY\right\} +\frac{p_{n}}{n}.
\end{equation*}%
Note $X^{2},Y^{2}\sim \chi ^{2}(1)$ and $\rho _{n}^{-n}=o\left(
c_{n}^{-1}\right) $. If $\lim\limits_{n\rightarrow \infty }\dfrac{p_{n}}{%
\rho _{n}^{2n}}=\pi $,

\begin{equation*}
\dfrac{n}{k_{n}}\left( IC_{1}-IC_{0}\right) \Rightarrow 
\begin{cases}
\pi -\dfrac{1}{4}\chi ^{2}(1), & \text{if }\pi \in \lbrack 0,\infty ) \\ 
1, & \text{if }\pi =\infty%
\end{cases}%
,
\end{equation*}%
where%
\begin{equation*}
k_{n}=%
\begin{cases}
\rho _{n}^{2n}, & \text{if }\pi \in \lbrack 0,\infty ) \\ 
p_{n}, & \text{if }\pi =\infty%
\end{cases}%
.
\end{equation*}

\subsection{Proof of Theorem \protect\ref{thm02c}}

When the true DGP is EX model, we have 
\begin{align*}
IC_{0}& =\log \widehat{\sigma }_{0}^{2}=\log \left\{ \dfrac{1}{n}%
\sum_{t=1}^{n}\left( X_{t}-X_{t-1}\right) ^{2}\right\} \\
& =\log \left\{ \dfrac{1}{n}\sum_{t=1}^{n}\left[ \left( \rho -1\right)
X_{t-1}+u_{t}\right] ^{2}\right\} \\
& =\log \left\{ \dfrac{1}{n}\left( \rho -1\right)
^{2}\sum_{t=1}^{n}X_{t-1}^{2}+\dfrac{2}{n}\left( \rho -1\right)
\sum_{t=1}^{n}X_{t-1}u_{t}+\dfrac{1}{n}\sum_{t=1}^{n}u_{t}^{2}\right\} .
\end{align*}

By results established in Anderson (1959), we know 
\begin{align}
\dfrac{1}{\rho ^{2n}}\sum_{t=1}^{n}X_{t-1}^{2}& \Rightarrow \dfrac{\sigma
^{2}Y^{2}}{\left( \rho ^{2}-1\right) ^{2}},  \label{eqES01} \\
\dfrac{1}{\rho ^{n}}\sum_{t=1}^{n}X_{t-1}u_{t}& \Rightarrow \dfrac{\sigma
^{2}XY}{\rho ^{2}-1},  \label{eqES02} \\
\dfrac{\rho ^{n}}{\rho ^{2}-1}\left( \widehat{\rho }-\rho \right) &
\Rightarrow \mathcal{C},  \label{eqES03}
\end{align}%
where $X,Y\overset{iid}{\sim }\mathcal{N}\left( 0,1\right) $ and $\mathcal{C}
$ is a standard Cauchy variate. Then we have 
\begin{align}
IC_{0}& =\log \left\{ \dfrac{\sigma ^{2}\rho ^{2n}}{n\left( \rho +1\right)
^{2}}X^{2}+\dfrac{2\sigma ^{2}\rho ^{n}}{n\left( \rho +1\right) }XY+\sigma
^{2}\right\}  \notag \\
& =\log \sigma ^{2}+\log \left\{ \dfrac{\rho ^{2n}}{n\left( \rho +1\right)
^{2}}X^{2}+\dfrac{2\rho ^{n}}{n\left( \rho +1\right) }XY+1\right\} .
\label{eq09}
\end{align}

For the OLS estimator for the general explosive series, we have 
\begin{align*}
IC_{1}& =\log \widehat{\sigma }_{1}^{2}+\frac{p_{n}}{n} \\
& =\log \left\{ n^{-1}\sum_{t=1}^{n}\left( X_{t}-\widehat{\rho }%
X_{t-1}\right) ^{2}\right\} +\frac{p_{n}}{n} \\
& =\log \left\{ \frac{1}{n}\sum_{t=1}^{n}\left[ \left( \rho -\widehat{\rho }%
\right) X_{t-1}+u_{t}\right] ^{2}\right\} +\frac{p_{n}}{n} \\
& =\log \left\{ \frac{1}{n}\left( \rho -\widehat{\rho }\right)
^{2}\sum_{t=1}^{n}X_{t-1}^{2}+\frac{2}{n}\left( \rho -\widehat{\rho }\right)
\sum_{t=1}^{n}X_{t-1}u_{t}+\frac{1}{n}\sum_{t=1}^{n}u_{t}^{2}\right\} +\frac{%
p_{n}}{n}.
\end{align*}

By equation (\ref{eqES01}) to (\ref{eqES03}), we have 
\begin{equation}
IC_{1}=\log \sigma ^{2}+\log \left\{ 1-\frac{1}{n}X^{2}\right\} +\frac{p_{n}%
}{n}.  \label{eq10}
\end{equation}

Now, by equation (\ref{eq09}) and (\ref{eq10}), we obtain 
\begin{equation*}
IC_{1}-IC_{0}=\log \left\{ 1-\frac{1}{n}X^{2}\right\} -\log \left\{ 1+\dfrac{%
2\rho ^{n}}{n\left( \rho +1\right) }XY+\dfrac{\rho ^{2n}}{n\left( \rho
+1\right) ^{2}}X^{2}\right\} +\frac{p_{n}}{n}.
\end{equation*}%
Since $\lim\limits_{n\rightarrow \infty }\dfrac{p_{n}}{\rho ^{2n}}=\pi $, we
have 
\begin{equation*}
\dfrac{n}{k_{n}}\left( IC_{1}-IC_{0}\right) \Rightarrow 
\begin{cases}
\pi -\dfrac{1}{(1+\rho )^{2}}\chi ^{2}(1), & \text{if }\pi \in \lbrack
0,\infty ) \\ 
1, & \text{if }\pi =\infty%
\end{cases}%
,
\end{equation*}%
where 
\begin{equation*}
k_{n}=%
\begin{cases}
\rho ^{2n}, & \text{if }\pi \in \lbrack 0,\infty ) \\ 
p_{n}, & \text{if }\pi =\infty%
\end{cases}%
.
\end{equation*}

\subsection{Proof of Proposition \protect\ref{prop05}}

When the true DGP is ME model, we have $0<c<\infty $, and%
\begin{align*}
IC_{0}& =\log \widehat{\sigma }_{0}^{2}=\log \left\{ \frac{1}{n}%
\sum_{t=1}^{n}\left( X_{t}-X_{t-1}\right) ^{2}\right\} \\
& =\log \left\{ \frac{1}{n}\sum_{t=1}^{n}\left[ \left( \rho _{n}-1\right)
X_{t-1}+u_{t}\right] ^{2}\right\} \\
& =\log \left\{ \frac{1}{n}\left( \rho _{n}-1\right)
^{2}\sum_{t=1}^{n}X_{t-1}^{2}+\frac{2}{n}\left( \rho _{n}-1\right)
\sum_{t=1}^{n}X_{t-1}u_{t}+\frac{1}{n}\sum_{t=1}^{n}u_{t}^{2}\right\} .
\end{align*}

When the process is initialized at $X_{0}$, by Lemma 5 in Magdalinos (2012),
we know 
\begin{equation}
\frac{c_{n}^{2}\rho _{n}^{-2n}}{\omega ^{2}n^{2}}\sum_{t=1}^{n}X_{t-1}^{2}%
\Rightarrow \frac{\sigma ^{2}}{4}Z^{2},  \label{eqMED01}
\end{equation}%
and 
\begin{equation}
\frac{c_{n}\rho _{n}^{-n}}{\omega ^{2}n}\sum_{t=1}^{n}X_{t-1}u_{t}%
\Rightarrow \frac{\sigma ^{2}}{2}YZ,  \label{eqMED02}
\end{equation}%
where by Lemma 2 in Magdalinos (2012), we know $Y$ and $Z$ are independent $%
\mathcal{N}\left( 0,1\right) $ variates with $\omega ^{2}=\left(
\sum\nolimits_{j=0}^{\infty }F_{j}\right) ^{2}$.

Therefore, by Equation (\ref{eqMED01}) and (\ref{eqMED02}) we have 
\begin{align}
IC_{0}& \Rightarrow \log \left\{ \frac{\omega ^{2}\sigma ^{2}}{4n\rho
_{n}^{-2n}}Z^{2}+\frac{\omega ^{2}\sigma ^{2}}{n\rho _{n}^{-n}}YZ+\sigma
^{2}+o_{p}(n^{-1})\right\}  \notag \\
& =\log \sigma ^{2}+\log \left\{ 1+\frac{\omega ^{2}}{n\rho _{n}^{-n}}YZ+%
\frac{\omega ^{2}}{4n\rho _{n}^{-2n}}Z^{2}+o_{p}(n^{-1})\right\} .
\label{eq12}
\end{align}

We also know from Magdalinos (2012) that 
\begin{equation}
\frac{n\rho _{n}^{n}}{2c_{n}}\left( \widehat{\rho }_{n}-\rho _{n}\right)
\Rightarrow \mathcal{C}.  \label{eqMED03}
\end{equation}%
Hence, 
\begin{align}
IC_{1}& =\log \widehat{\sigma }_{1}^{2}+\frac{p_{n}}{n}  \notag \\
& =\log \left\{ n^{-1}\sum_{t=1}^{n}\left( X_{t}-\widehat{\rho }%
_{n}X_{t-1}\right) ^{2}\right\} +\frac{p_{n}}{n}  \notag \\
& =\log \left\{ \frac{1}{n}\sum_{t=1}^{n}\left[ \left( \rho _{n}-\widehat{%
\rho }_{n}\right) X_{t-1}+u_{t}\right] ^{2}\right\} +\frac{p_{n}}{n}  \notag
\\
& =\log \left\{ \frac{1}{n}\left( \rho _{n}-\widehat{\rho }_{n}\right)
^{2}\sum_{t=1}^{n}X_{t-1}^{2}+\frac{2}{n}\left( \rho _{n}-\widehat{\rho }%
_{n}\right) \sum_{t=1}^{n}X_{t-1}u_{t}+\frac{1}{n}\sum_{t=1}^{n}u_{t}^{2}%
\right\} +\frac{p_{n}}{n}  \notag \\
& \Rightarrow \log \left\{ -\dfrac{\omega ^{2}\sigma ^{2}}{n}Y^{2}+\sigma
^{2}\right\} +\frac{p_{n}}{n}  \notag \\
& =\log \sigma ^{2}+\log \left\{ 1-\dfrac{\omega ^{2}}{n}Y^{2}\right\} +%
\frac{p_{n}}{n}.  \label{eq13}
\end{align}

Therefore, by Equation (\ref{eq12}) and (\ref{eq13}), we have 
\begin{equation*}
IC_{1}-IC_{0}\Rightarrow \log \left\{ 1-\dfrac{\omega ^{2}}{n}Y^{2}\right\}
-\log \left\{ 1+\frac{\omega ^{2}}{n\rho _{n}^{-n}}YZ+\frac{\omega ^{2}}{%
4n\rho _{n}^{-2n}}Z^{2}\right\} +\frac{p_{n}}{n}.
\end{equation*}

Note $Y^{2},Z^{2}\sim \chi ^{2}(1)$ and $\rho _{n}^{-n}=o\left(
c_{n}^{-1}\right) $. If $\lim\limits_{n\rightarrow \infty }\dfrac{p_{n}}{%
\rho _{n}^{2n}}=\pi $,

\begin{equation*}
\dfrac{n}{k_{n}}\left( IC_{1}-IC_{0}\right) \Rightarrow 
\begin{cases}
\pi -\dfrac{\omega ^{2}}{4}\chi ^{2}(1), & \text{if }\pi \in \lbrack
0,\infty ) \\ 
1, & \text{if }\pi =\infty%
\end{cases}%
,
\end{equation*}%
where%
\begin{equation*}
k_{n}=%
\begin{cases}
\rho _{n}^{2n}, & \text{if }\pi \in \lbrack 0,\infty ) \\ 
p_{n}, & \text{if }\pi =\infty%
\end{cases}%
.
\end{equation*}

\subsection{Proof of Theorem \protect\ref{thm03}}

When the true DGP is the UR\ model, we have 
\begin{equation*}
IC_{0}=\log \breve{\sigma}_{0}^{2}=\log \left\{ \frac{1}{n}%
\sum_{t=1}^{n}u_{t}^{2}\right\} =\log \sigma ^{2}.
\end{equation*}%
Also, we have 
\begin{align*}
IC_{1}& =\log \breve{\sigma}_{1}^{2}+\dfrac{p_{n}}{n}=\log \left\{ \frac{1}{n%
}\sum_{t=1}^{n}\left( X_{t}-\breve{\rho}X_{t-1}\right) ^{2}\right\} +\dfrac{%
p_{n}}{n} \\
& =\log \left\{ \frac{1}{n}\sum_{t=1}^{n}\left[ \left( 1-\breve{\rho}\right)
X_{t-1}+u_{t}\right] ^{2}\right\} +\dfrac{p_{n}}{n} \\
& =\log \left\{ \frac{1}{n}\left( 1-\breve{\rho}\right)
^{2}\sum_{t=1}^{n}X_{t-1}^{2}+\frac{2}{n}\left( 1-\breve{\rho}\right)
\sum_{t=1}^{n}X_{t-1}u_{t}+\frac{1}{n}\sum_{t=1}^{n}u_{t}^{2}\right\} +%
\dfrac{p_{n}}{n}.
\end{align*}

According to Phillips (2012), we have 
\begin{equation*}
\breve{\rho}-1\Rightarrow 
\begin{cases}
\dfrac{1}{n}h^{-1}\left( \dfrac{\int_{0}^{1}BdB}{\int_{0}^{1}B^{2}}\right) ,
& \text{if }\tau =0 \\ 
\dfrac{1}{n}h^{-1}\left( \dfrac{\int_{0}^{1}B_{\tau }dB}{\int_{0}^{1}B_{\tau
}^{2}}\right) , & \text{if }\tau \in (0,\infty ) \\ 
\dfrac{1}{\sqrt{n^{2}/c_{n}}}h^{-1}\left( \mathcal{C}\right) , & \text{if }%
\tau =\infty%
\end{cases}%
,
\end{equation*}%
where $h(c)$ was defined in Section 4.

According to Phillips and Magdalinos (2009), we have 
\begin{align*}
\frac{1}{n^{2}}\sum_{t=1}^{n}X_{t-1}^{2}& \Rightarrow 
\begin{cases}
\sigma ^{2}\int_{0}^{1}B^{2}, & \text{if }\tau =0 \\ 
\sigma ^{2}\int_{0}^{1}B_{\tau }^{2}, & \text{if }\tau \in (0,\infty ) \\ 
\sigma ^{2}B_{0}(1)^{2}/c_{n}, & \text{if }\tau =\infty%
\end{cases}%
, \\
\frac{1}{n}\sum_{t=1}^{n}X_{t-1}u_{t}& \Rightarrow 
\begin{cases}
\sigma ^{2}\int_{0}^{1}BdB, & \text{if }\tau =0 \\ 
\sigma ^{2}\int_{0}^{1}B_{\tau }dB, & \text{if }\tau \in (0,\infty ) \\ 
\sqrt{1/c_{n}}\sigma ^{2}B(1)B_{0}(1), & \text{if }\tau =\infty%
\end{cases}%
.
\end{align*}%
Therefore, we have 
\begin{equation*}
IC_{0}-IC_{1}\Rightarrow 
\begin{cases}
-\log \left\{ \dfrac{\int_{0}^{1}B^{2}}{n}h^{-1}\left( \dfrac{\int_{0}^{1}BdB%
}{\int_{0}^{1}B^{2}}\right) ^{2}-\dfrac{2\left( \int_{0}^{1}BdB\right) }{n}%
h^{-1}\left( \dfrac{\int_{0}^{1}BdB}{\int_{0}^{1}B^{2}}\right) +1\right\} -%
\dfrac{p_{n}}{n} \\ 
-\log \left\{ \dfrac{\int_{0}^{1}B_{\tau }^{2}}{n}h^{-1}\left( \dfrac{%
\int_{0}^{1}B_{\tau }dB}{\int_{0}^{1}B_{\tau }^{2}}\right) ^{2}-\dfrac{%
2\left( \int_{0}^{1}B_{\tau }dB\right) }{n}h^{-1}\left( \dfrac{%
\int_{0}^{1}B_{\tau }dB}{\int_{0}^{1}B_{\tau }^{2}}\right) +1\right\} -%
\dfrac{p_{n}}{n} \\ 
-\log \left\{ \dfrac{1}{n}h^{-1}\left( \mathcal{C}\right) ^{2}B_{0}(1)^{2}-%
\dfrac{2}{n}h^{-1}\left( \mathcal{C}\right) B(1)B_{0}(1)+1\right\} -\dfrac{%
p_{n}}{n}%
\end{cases}%
.
\end{equation*}

\subsection{Proof of Theorem \protect\ref{thm04a}}

When the true DGP is the LTUE model, we have $0<c<\infty $. There is no
difference between $IC_{0}$ based on the OLS estimator and that based on the
indirect inference estimator. For $IC_{1}$, we have 
\begin{align*}
IC_{1}& =\log \breve{\sigma}_{1}^{2}+\frac{p_{n}}{n} \\
& =\log \left\{ n^{-1}\sum_{t=1}^{n}\left( X_{t}-\breve{\rho}%
_{n}X_{t-1}\right) ^{2}\right\} +\frac{p_{n}}{n} \\
& =\log \left\{ \frac{1}{n}\sum_{t=1}^{n}\left[ \left( \rho _{n}-\breve{\rho}%
_{n}\right) X_{t-1}+u_{t}\right] ^{2}\right\} +\frac{p_{n}}{n} \\
& =\log \left\{ \frac{1}{n}\left( \rho _{n}-\breve{\rho}_{n}\right)
^{2}\sum_{t=1}^{n}X_{t-1}^{2}+\frac{2}{n}\left( \rho _{n}-\breve{\rho}%
_{n}\right) \sum_{t=1}^{n}X_{t-1}u_{t}+\frac{1}{n}\sum_{t=1}^{n}u_{t}^{2}%
\right\} +\frac{p_{n}}{n}.
\end{align*}

By the limit theory for the indirect inference estimator developed in
Phillips (2012), we have 
\begin{equation}
n\left( \breve{\rho}_{n}-\rho _{n}\right) \Rightarrow h^{-1}\left( \dfrac{%
\int_{0}^{1}J_{c}dB}{\int_{0}^{1}J_{c}^{2}}+c\right) -c.  \label{eq05}
\end{equation}

By equation (\ref{eqLU01}), (\ref{eqLU02}) and (\ref{eq05}), we have 
\begin{alignat}{2}
IC_{1}& \Rightarrow & & \log \left\{ 1-\dfrac{2}{n}\left[ h^{-1}\left( 
\dfrac{\int_{0}^{1}J_{c}dB}{\int_{0}^{1}J_{c}^{2}}+c\right) -c\right]
\int_{0}^{1}J_{c}dB+\dfrac{1}{n}\left[ h^{-1}\left( \dfrac{%
\int_{0}^{1}J_{c}dB}{\int_{0}^{1}J_{c}^{2}}+c\right) -c\right]
^{2}\int_{0}^{1}J_{c}^{2}\right\}  \notag \\
& {} & & +\log \sigma ^{2}+\frac{p_{n}}{n}.  \label{eq06}
\end{alignat}

Therefore, by equation (\ref{eq03}) and (\ref{eq05}), we have 
\begin{alignat*}{2}
IC_{1}-IC_{0}& \Rightarrow & & \log \left\{ 1-\dfrac{2\int_{0}^{1}J_{c}dB}{n}%
\left[ h^{-1}\left( \dfrac{\int_{0}^{1}J_{c}dB}{\int_{0}^{1}J_{c}^{2}}%
+c\right) -c\right] +\dfrac{\int_{0}^{1}J_{c}^{2}}{n}\left[ h^{-1}\left( 
\dfrac{\int_{0}^{1}J_{c}dB}{\int_{0}^{1}J_{c}^{2}}+c\right) -c\right]
^{2}\right\} \\
& {} & & -\log \left\{ 1+\dfrac{2c}{n}\int_{0}^{1}J_{c}dB+\dfrac{c^{2}}{n}%
\int_{0}^{1}J_{c}^{2}\right\} +\dfrac{p_{n}}{n}.
\end{alignat*}

When $p_{n}=2$, as $n\rightarrow \infty $ we have 
\begin{equation*}
n\left( IC_{1}-IC_{0}\right) \Rightarrow 2-\vartheta ^{2}.
\end{equation*}%
where 
\begin{equation*}
\vartheta ^{2}\equiv 2h^{-1}\left( \dfrac{\int_{0}^{1}J_{c}dB}{%
\int_{0}^{1}J_{c}^{2}}+c\right) \left(
\int_{0}^{1}J_{c}dB+c\int_{0}^{1}J_{c}^{2}\right) -h^{-1}\left( \dfrac{%
\int_{0}^{1}J_{c}dB}{\int_{0}^{1}J_{c}^{2}}+c\right)
^{2}\int_{0}^{1}J_{c}^{2}.
\end{equation*}%
When $p_{n}\rightarrow \infty $ and $\dfrac{p_{n}}{n}\rightarrow 0,$ we have 
\begin{equation*}
\dfrac{n}{p_{n}}\left( IC_{1}-IC_{0}\right) \Rightarrow 1.
\end{equation*}

\subsection{Proof of Theorem \protect\ref{thm04b}}

When the true DGP is the ME model, we have 
\begin{align*}
IC_{0}& =\log \breve{\sigma}_{0}^{2}=\log \left\{ \frac{1}{n}%
\sum_{t=1}^{n}\left( X_{t}-X_{t-1}\right) ^{2}\right\} \\
& =\log \left\{ \frac{1}{n}\sum_{t=1}^{n}\left[ \left( \rho _{n}-1\right)
X_{t-1}+u_{t}\right] ^{2}\right\} \\
& =\log \left\{ \frac{1}{n}\left( \rho _{n}-1\right)
^{2}\sum_{t=1}^{n}X_{t-1}^{2}+\frac{2}{n}\left( \rho _{n}-1\right)
\sum_{t=1}^{n}X_{t-1}u_{t}+\frac{1}{n}\sum_{t=1}^{n}u_{t}^{2}\right\} .
\end{align*}

By equation (\ref{eqPM01}) and (\ref{eqPM02}) we have 
\begin{align}
IC_{0}& \Rightarrow \log \left\{ \frac{\sigma ^{2}}{4n\rho _{n}^{-2n}}Y^{2}+%
\frac{\sigma ^{2}}{n\rho _{n}^{-n}}XY+\sigma ^{2}\right\}  \notag \\
& =\log \sigma ^{2}+\log \left\{ \frac{1}{4n\rho _{n}^{-2n}}Y^{2}+\frac{1}{%
n\rho _{n}^{-n}}XY+1\right\} .  \label{eq07}
\end{align}

Similarly, for $IC_{1}$ based on the indirect inference estimator, we have 
\begin{align*}
IC_{1}& =\log \breve{\sigma}_{1}^{2}+\frac{p_{n}}{n} \\
& =\log \left\{ n^{-1}\sum_{t=1}^{n}\left( X_{t}-\breve{\rho}%
_{n}X_{t-1}\right) ^{2}\right\} +\frac{p_{n}}{n} \\
& =\log \left\{ \frac{1}{n}\sum_{t=1}^{n}\left[ \left( \rho _{n}-\breve{\rho}%
_{n}\right) X_{t-1}+u_{t}\right] ^{2}\right\} +\frac{p_{n}}{n} \\
& =\log \left\{ \frac{1}{n}\left( \rho _{n}-\breve{\rho}_{n}\right)
^{2}\sum_{t=1}^{n}X_{t-1}^{2}+\frac{2}{n}\left( \rho _{n}-\breve{\rho}%
_{n}\right) \sum_{t=1}^{n}X_{t-1}u_{t}+\frac{1}{n}\sum_{t=1}^{n}u_{t}^{2}%
\right\} +\frac{p_{n}}{n}
\end{align*}

Using the results in Phillips (2012) , equation (\ref{eqPM01}) and (\ref%
{eqPM02}), we obtain 
\begin{align}
IC_{1}& \Rightarrow \log \left\{ \dfrac{1}{n}\dfrac{4c_{n}^{2}}{n^{2}\rho
_{n}^{2n}}\left( \mathcal{C}+O\left( \dfrac{1}{2c_{n}}\right) \right) ^{2}%
\dfrac{n^{2}\sigma ^{2}}{4c_{n}^{2}\rho _{n}^{-2n}}Y^{2}-\dfrac{2}{n}\dfrac{%
2c_{n}}{n\rho _{n}^{n}}\left( \mathcal{C}+O\left( \dfrac{1}{2c_{n}}\right)
\right) \dfrac{n\sigma ^{2}}{2c_{n}\rho _{n}^{-n}}XY+\sigma ^{2}\right\} +%
\dfrac{p_{n}}{n}  \notag \\
& \Rightarrow \log \sigma ^{2}+\log \left\{ 1-\dfrac{1}{n}X^{2}+O\left( 
\dfrac{1}{c_{n}n}\right) \right\} +\dfrac{p_{n}}{n}.  \label{eq08}
\end{align}

Therefore, the similar results to those in Theorem \ref{thm02b} are obtained.

\subsection{Proof of Theorem \protect\ref{thm04c}}

When the true DGP is the EX model, for the indirect inference estimator, we
know that for $IC_{0}$, it is the same as OLS estimator. Therefore, we only
need to derive the $IC_{1}$. Note that for $IC_{1}$, we have 
\begin{equation*}
IC_{1}=\log \left\{ \frac{1}{n}\left( \rho -\breve{\rho}\right)
^{2}\sum_{t=1}^{n}X_{t-1}^{2}+\frac{2}{n}\left( \rho -\breve{\rho}\right)
\sum_{t=1}^{n}X_{t-1}u_{t}+\frac{1}{n}\sum_{t=1}^{n}u_{t}^{2}\right\} +\frac{%
p_{n}}{n}.
\end{equation*}

According to the results in Phillips (2012), for $|\rho |>1$, we know the
binding function for $\rho $ is 
\begin{equation*}
b_{n}(\rho )=\rho +O(\rho ^{-n}).
\end{equation*}

Therefore, we obtain 
\begin{align}
IC_{1}& =\log \sigma ^{2}+\log \left\{ \frac{1}{n}\left( \mathcal{C}+O\left( 
\dfrac{1}{\rho ^{2}-1}\right) \right) ^{2}Y^{2}-\frac{2}{n}\left( \mathcal{C}%
+O\left( \dfrac{1}{\rho ^{2}-1}\right) \right) XY+1\right\} +\frac{p_{n}}{n}
\notag \\
& =\log \sigma ^{2}+\log \left\{ 1-\dfrac{1}{n}X^{2}+O\left( \dfrac{1}{%
n\left( \rho ^{2}-1\right) }\right) \right\} +\frac{p_{n}}{n}.  \label{eq11}
\end{align}

Now, by equation (\ref{eq09}) and (\ref{eq11}), we obtain 
\begin{equation*}
IC_{1}-IC_{0}=\log \left\{ 1-\frac{1}{n}X^{2}+O\left( \dfrac{1}{n\left( \rho
^{2}-1\right) }\right) \right\} -\log \left\{ 1+\dfrac{2\rho ^{n}}{n\left(
\rho +1\right) }XY+\dfrac{\rho ^{2n}}{n\left( \rho +1\right) ^{2}}%
X^{2}\right\} +\frac{p_{n}}{n}.
\end{equation*}

Since $\lim\limits_{n\rightarrow \infty }\dfrac{p_{n}}{\rho ^{2n}}=\pi $, we
have 
\begin{equation*}
\dfrac{n}{k_{n}}\left( IC_{1}-IC_{0}\right) \Rightarrow 
\begin{cases}
\pi -\dfrac{1}{(1+\rho )^{2}}\chi ^{2}(1), & \text{if }\pi \in \lbrack
0,\infty ) \\ 
1, & \text{if }\pi =\infty%
\end{cases}%
,
\end{equation*}%
where 
\begin{equation*}
k_{n}=%
\begin{cases}
\rho ^{2n}, & \text{if }\pi \in \lbrack 0,\infty ) \\ 
p_{n}, & \text{if }\pi =\infty%
\end{cases}%
.
\end{equation*}

\end{document}